\documentclass[12pt,a4wide]{article}
\usepackage{amsfonts}
\usepackage{amssymb}
\usepackage{amsthm}
\usepackage{amsmath}
\usepackage[all]{xy}
\usepackage{enumerate}
\usepackage{fancyhdr}
\usepackage{a4wide}

\newtheorem{theorem}{Theorem}[section]
\newtheorem{proposition}[theorem]{Proposition}
\newtheorem{corollary}[theorem]{Corollary}
\newtheorem{lemma}[theorem]{Lemma}
\theoremstyle{definition}
\newtheorem*{notation}{Notation}
\newtheorem*{Beweis}{Proof}
\newtheorem{definition}[theorem]{Definition}
\newtheorem{punto}[theorem]{}
\theoremstyle{remark}
\newtheorem{remark}[theorem]{Remark}
\newtheorem{ex}[theorem]{Example}
\newtheorem{exs}[theorem]{Examples}
\newtheorem{remarks}[theorem]{Remarks}
\newcommand{\V}{\mathcal V}
\newcommand{\W}{\mathcal W}
\input{tcilatex}
\begin{document}

\title{\textbf{Semiunital Semimonoidal Categories}\\
(\textbf{Applications to Semirings and Semicorings)}}
\author{{\textbf{Jawad Abuhlail}}\thanks{%
The author would like to acknowledge the support provided by the Deanship of
Scientific Research (DSR) at King Fahd University of Petroleum $\&$ Minerals
(KFUPM) for funding this work through project No. IN100008.} \\
Department of Mathematics and Statistics\\
Box 5046, KFUPM, 31261 Dhahran, KSA\\
abuhlail@kfupm.edu.sa}
\date{\today}
\maketitle

\begin{abstract}
The category $_{A}\mathbb{S}_{A}$ of bisemimodules over a semialgebra $A,$
with the so called \emph{Takahashi's tensor product} $-\boxtimes _{A}-,$ is
semimonoidal but \emph{not} monoidal. Although not a unit in $_{A}\mathbb{S}%
_{A},$ the base semialgebra $A$ has properties of a \emph{semiunit} (in a
sense which we clarify in this note). Motivated by this interesting example,
we investigate \emph{semiunital semimonoidal categories }$(\mathcal{V}%
,\bullet ,\mathbf{I})$\emph{\ }as a framework for studying notions like
\emph{semimonoids} (\emph{semicomonoids}) as well as a notion of monads
(comonads) which we call $\mathbb{J}$-\emph{monads} ($\mathbb{J}$-\emph{%
comonads}) with respect to the endo-functor $\mathbb{J}:=\mathbf{I}\bullet
-\simeq -\bullet \mathbf{I}:\mathcal{V}\longrightarrow \mathcal{V}.$ This
motivated also introducing a more generalized notion of monads (comonads)
in arbitrary categories with respect to arbitrary endo-functors.
Applications to the semiunital semimonoidal variety $(_{A}\mathbb{S}%
_{A},\boxtimes _{A},A)$ provide us with examples of semiunital $A$-semirings
(semicounital $A$-semicorings) and semiunitary semimodules (semicounitary
semicomodules) which extend the classical notions of unital rings (counital
corings) and unitary modules (counitary comodules).
\end{abstract}

\section{Introduction}

\qquad A \emph{semiring} is, roughly speaking, a ring not necessarily with
subtraction. The first natural example of a semiring is the set $\mathbb{N}%
_{0}$ of non-negative integers. Other examples include the set $\mathrm{Ideal%
}(R)$ of (two-sided) ideals of any associative ring $R$ and distributive
complete lattices. A \emph{semimodule} is, roughly speaking, a module not
necessarily with subtraction. The category of Abelian groups is nothing but
the category of modules over $\mathbb{Z};$ similarly, the category of
commutative monoids is nothing but the category of semimodules over $\mathbb{%
N}_{0}.$

Semirings were studied by many algebraists beginning with Dedekind \cite%
{Ded1894}. Since the sixties of the last century, they were shown to have
significant applications in several areas as Automata Theory, Optimization
Theory, Tropical Geometry and Idempotent Analysis (for more, see \cite%
{Gol1999a}). Recently, Durov \cite{Dur2007} demonstrated that semirings are
in one-to-one correspondence with the \emph{algebraic additive monads} on
the category $\mathbf{Set}$ of sets. The theory of semimodules over
semirings was developed by many authors including Takahashi, Patchkoria and
Katsov (\emph{e.g.} \cite{Tak1981}, \cite{Tak1982a}, \cite{Pat2006}, \cite%
{Kat1997}).

A strong connection between corings \cite{Swe1975} over a ring $A$
(coalgebras in the monoidal category $_{A}\mathbf{Mod}_{A}$ of bimodules
over $A$) and their comodules on one side and comonads induced by the tensor
product $-\otimes _{A}-$ and their comodules on the other side has been
realized by several authors (\emph{e.g.} \cite{BW2003}). Moreover, the
theory of monads and comonads in (autonomous) monoidal categories received
increasing attention in the last decade and extensions to arbitrary
categories were carried out in several recent papers (\emph{e.g.} \cite%
{BBW2009}).

Using the so called Takahashi's \emph{tensor-like product} $-\boxtimes _{A}-$
of semimodules over an associative semiring $A$ \cite{Tak1982a}, notions of
\emph{semiunital semirings }and \emph{semicounital semicorings} were
introduced by the author in 2008. However, these could not be realized as
monoids (comonoids) in the category $_{A}\mathbb{S}_{A}$ of $(A,A)$%
-bisemimodules. This is mainly due to the fact that the category $(_{A}%
\mathbb{S}_{A},\boxtimes _{A},A)$ is \emph{not} monoidal in general (an
alternative tensor product $-\otimes _{A}-$ was recalled by Katsov in \cite%
{Kat1997}; in fact $(_{A}\mathbb{S}_{A},\otimes _{A},A)$ is monoidal. For
the relation between $-\otimes _{A}-$ and $-\boxtimes _{A}-,$ see \cite{Abu}%
). Motivated by the desire to fix this defect, we introduce and investigate
a notion of \emph{semiunital semimonoidal categories} with prototype $(_{A}%
\mathbb{S}_{A},\boxtimes ,A)$ and investigate \emph{semimonoids} (\emph{%
semicomonoids}) in such categories as well as their categories of \emph{%
semimodules} (\emph{semicomodules}). In particular, we realize our
semiunital $A$-semirings (semicounital $A$-semicorings) as semimonoids
(semicomonoids) in $(_{A}\mathbb{S}_{A},\boxtimes ,A).$ Moreover, we
introduce and study $\mathbb{J}$-\emph{monads }($\mathbb{J}$-\emph{comonads})%
\emph{\ }in any arbitrary category $\mathfrak{A},$ where $\mathbb{J}:%
\mathfrak{A}\longrightarrow \mathfrak{A}$ is an endo-functor, and apply them
to semiunital semimonoidal categories in general and to $_{A}\mathbb{S}_{A}$
in particular. Our results extend recent ones on monoids (comonoids) in
monoidal categories as well as monads (comonads) in arbitrary categories to
semimonoids (semicomonoids) in semiunital semimonoidal categories as well as
$\mathbb{J}$-monads ($\mathbb{J}$-comonads) in arbitrary categories.

Throughout, $\mathbb{I}$ denotes the identity endo-functor on the category
under consideration. The paper is organized as follows. After this
introduction, we present in Section 2 our (generalized) notion of $\mathbb{J}
$-monads and $\mathbb{J}$-comonads in arbitrary categories. In Section 3, we
introduce and investigate \emph{semiunits} in semimonoidal categories. In
Section 4, we introduce semimonoids (semicomonoids) in semiunital
semimonoidal categories as well as their categories of semimodules
(semicomodules). Moreover, we prove two \emph{reconstruction} results,
namely Theorems \ref{thm-monad} and \ref{thm-comonad}. In Section 5, we
consider the semiunital semimonoidal category (variety) of bisemimodules $%
_{A}\mathbb{S}_{A}$ over a semialgebra $A$ which provides us with a rich
source of concrete examples for applying our results. As mentioned above,
these concrete examples were the main motivation behind introducing all the
abstract notions in this paper. Further investigations of $\mathbb{J}$\emph{%
-bimonads} and \emph{Hopf} $\mathbb{J}$\emph{-monads} as well as \emph{%
bisemimonoids} and \emph{Hopf semimonoids} in semiunital semimonoidal
categories will be the subject of a forthcoming paper.

\section{Monads and Comonads}

\qquad Recall first the so called \emph{Godement product}\ of natural
transformations between functors:

\begin{punto}
Let $\mathfrak{A},\mathfrak{B},\mathfrak{C}$ be any categories. Any natural
transformations $\mathbf{\psi }:F\longrightarrow G$ and $\mathbf{\phi }%
:F^{\prime }\longrightarrow G^{\prime }$ of functors $\mathfrak{A}\overset{%
F,G}{\longrightarrow }\mathfrak{B}\overset{F^{\prime },G^{\prime }}{%
\longrightarrow }\mathfrak{C}$ can be multiplied using the \emph{Godement
product} to yield a natural transformation $\mathbf{\phi \psi }:F^{\prime
}F\longrightarrow G^{\prime }G,$ where
\begin{equation}
\mathbf{\phi }_{G(X)}\circ F^{\prime }(\mathbf{\psi }_{X})=(\mathbf{\phi
\psi })_{X}=G^{\prime }(\mathbf{\psi }_{X})\circ \mathbf{\phi }_{F(X)}\text{
for every }X\in \mathfrak{A}.  \label{Godement}
\end{equation}

Moreover, if $\mathfrak{A}\overset{H}{\longrightarrow }\mathfrak{B}\overset{%
H^{\prime }}{\longrightarrow }\mathfrak{C}$ are functors and $\mathbf{\delta
}:G\longrightarrow H,$ $\mathbf{\theta }:G^{\prime }\longrightarrow
H^{\prime }$ are natural transformations, then the following \emph{%
interchange law }holds%
\begin{equation}
(\mathbf{\delta }\circ \mathbf{\psi })(\mathbf{\theta }\circ \mathbf{\phi }%
)=(\mathbf{\theta \delta })\circ (\mathbf{\phi \psi }).  \label{interchange}
\end{equation}
\end{punto}

\begin{punto}
Let $\mathfrak{A}$ and $\mathfrak{B}$ be categories, $L:\mathfrak{A}%
\longrightarrow \mathfrak{B},$ $R:\mathfrak{B}\longrightarrow \mathfrak{A}$
be functors and $\mathbb{J}:\mathfrak{A}\longrightarrow \mathfrak{A},$ $%
\mathbb{K}:\mathfrak{B}\longrightarrow \mathfrak{B}$ be endo-functors such
that $R\mathbb{K}\simeq \mathbb{J}R$ and $L\mathbb{J}\simeq \mathbb{K}L.$ We
say that $(L,R)$ is a\emph{\ }$(\mathbb{J},\mathbb{K})$\emph{-adjoint pair}
iff we have natural isomorphisms in $X\in \mathfrak{A}$ and $Y\in \mathfrak{B%
}:$%
\begin{equation*}
\mathfrak{B}(L\mathbb{J}(X),\mathbb{K}(Y))\simeq \mathfrak{A}(\mathbb{J}(X),R%
\mathbb{K}(Y)).
\end{equation*}%
For the special case $\mathbb{J}=\mathbb{I}_{\mathfrak{A}}$ and $\mathbb{K}=%
\mathbb{I}_{\mathfrak{B}},$ we recover the classical notion of adjoint pairs.
\end{punto}

Till the end of this section, $\mathfrak{A}$ is an arbitrary category.

\begin{punto}
Let $\mathbb{T}:\mathfrak{A}\longrightarrow \mathfrak{A}$ be an
endo-functor. An object $X\in \mathrm{Obj}(\mathfrak{A})$ is said to have a $%
\mathbb{T}$\emph{-action }or to be a $\mathbb{T}$\emph{-act} iff there is a
morphism $\varrho _{X}:\mathbb{T}(X)\longrightarrow X$ in $\mathfrak{A}.$
For two objects $X,X^{\prime }$ with $\mathbb{T}$-actions, we say that a
morphism $\varphi :X\longrightarrow X^{\prime }$ in $\mathfrak{A}$ is a
\emph{morphism of }$\mathbb{T}$\emph{-acts} iff the following diagram is
commutative%
\begin{equation*}
\xymatrix{ \mathbb{T}(X) \ar[d]_{\mathbb{T}(\varphi)} \ar[rr]^{\varrho_{X}}
& & X \ar[d]^{\varphi} \\ \mathbb{T}(X') \ar[rr]_{\varrho_{X'}} &&X'}
\end{equation*}%
The \emph{category of }$\mathbb{T}$\emph{-acts} is denoted by ${\mathbf{Act}}%
_{\mathbb{T}}.$ Dually, one can define the category $\mathbf{Coact}^{\mathbb{%
T}}$ of $\mathbb{T}$\emph{-coacts.}
\end{punto}

\begin{remark}
The objects of $\mathbf{Coact}^{\mathbb{F}},\ $where $\mathbb{F}:\mathbf{Set}%
\longrightarrow \mathbf{Set}$ is any endo-functor, play an important role in
logic and theoretical computer science. They are called $\mathbb{F}$\emph{%
-systems }(\emph{e.g.} \cite{Rut2000}). Some references call these $\mathbb{F%
}$\emph{-coalgebras }(\emph{e.g.} \cite{Gum1999}). For us, coalgebras are
always coassociative and counital unless something else is explicitly
specified.
\end{remark}

\subsection*{$\mathbb{J}$-Monads}

\begin{punto}
Let $\mathbb{J}:\mathfrak{A}\longrightarrow \mathfrak{A}$ be an
endo-functor. With a $\mathbb{J}$-\emph{monad} on $\mathfrak{A}$ we mean a
datum $(\mathbb{M},\mathbf{\mu },\mathbf{\omega },\mathbf{\nu };\mathbb{J})$
consisting of an endo-functor $\mathbb{M}:\mathfrak{A}\longrightarrow
\mathfrak{A}$ associated with natural transformations%
\begin{equation*}
\mathbf{\mu }:\mathbb{MM}\longrightarrow \mathbb{M},\text{ }\mathbf{\omega }:%
\mathbb{I}\longrightarrow \mathbb{J}\text{ and }\mathbf{\nu }:\mathbb{J}%
\longrightarrow \mathbb{M}
\end{equation*}%
such that the following diagrams are commutative
\begin{equation*}
\begin{tabular}{lll}
$\xymatrix{ \mathbb{MMM} \ar[dd]_{\mathbf{\mu} \mathbb{M}}
\ar[rr]^{\mathbb{M} \mathbf{\mu}} & & \mathbb{MM} \ar[dd]^{\mathbf{\mu}} \\
& & & \\ \mathbb{MM} \ar[rr]_{\mathbf{\mu}} \ar[rr]_{\mathbf{\mu}} & &
\mathbb{M}}$ & $\xymatrix{ \mathbb{MM} \ar[rr]^{\mathbf{\mu}} & &
\mathbb{IM} \ar[dd]^{\mathbf{\omega}\mathbb{M}} \\ & & & \\ \mathbb{J M}
\ar[uu]^{\mathbf{\nu}\mathbb{M}} & & \mathbb{JM} \ar@{=}[ll]}$ & $\xymatrix{
\mathbb{MM} \ar[rr]^{\mathbf{\mu}} & & \mathbb{M} \mathbb{I}
\ar[dd]^{\mathbb{M} \mathbf{\omega}} \\ & & & \\ \mathbb{M J}
\ar[uu]^{\mathbb{M}\mathbf{\nu}} & & \mathbb{M J} \ar@{=}[ll]}$%
\end{tabular}%
\end{equation*}%
\emph{i.e. }for every $X\in \mathfrak{A}$ we have%
\begin{equation*}
\mathbf{\mu }_{X}\circ \mathbb{M}(\mathbf{\mu }_{X})=\mathbf{\mu }_{X}\circ
\mathbf{\mu }_{\mathbb{M}(X)},\text{ }\mathbf{\nu }_{\mathbb{M}(X)}\circ
\mathbf{\omega }_{\mathbb{M}(X)}\circ \mathbf{\mu }_{X}=\mathbb{I}_{\mathbb{%
MM}(X)}\text{ and }\mathbb{M}(\mathbf{\nu }_{X})\circ \mathbb{M}(\mathbf{%
\omega }_{X})\circ \mathbf{\mu }_{X}=\mathbb{I}_{\mathbb{MM}(X)}.
\end{equation*}
\end{punto}

\begin{punto}
With $\mathbf{JMonad}_{\mathfrak{A}}$ we denote the category whose objects
are $\mathbb{J}$-monads, where $\mathbb{J}$ runs over the class of
endo-functors on $\mathfrak{A}.$ A morphism $(\mathbf{\varphi };\mathbf{\xi }%
):(\mathbb{M},\mathbf{\mu },\mathbf{\omega },\mathbf{\nu };\mathbb{J}%
)\longrightarrow (\mathbb{M}^{\prime },\mathbf{\mu }^{\prime },\mathbf{%
\omega }^{\prime },\mathbf{\nu }^{\prime };\mathbb{J}^{\prime })$ in this
category consists of natural transformations $\mathbf{\varphi }:\mathbb{M}%
\longrightarrow \mathbb{M}^{\prime }$ and $\mathbf{\xi }:\mathbb{J}%
\longrightarrow \mathbb{J}^{\prime }$ such that the following diagrams are
commutative%
\begin{equation*}
\begin{tabular}{lll}
$\xymatrix{ \mathbb{MM} \ar[d]_{\mathbf{\varphi \varphi}}
\ar[rr]^{\mathbf{\mu}} & & \mathbb{M} \ar[d]^{\mathbf{\varphi}} \\
{\mathbb{M' M'}} \ar[rr]_{{\mathbf{\mu}}'} & & \mathbb{M'} }$ &  & $%
\xymatrix{\mathbb{J} \ar[d]_{{\mathbf{\xi}}} \ar[rr]^{\mathbf{\nu}} & &
\mathbb{M} \ar[d]^{\mathbf{\varphi}}\\ \mathbb{J'} \ar[rr]_{{\mathbf{\nu}}'}
& & {\mathbb{M}'} }$%
\end{tabular}%
\end{equation*}%
\emph{i.e.} for every $X\in \mathfrak{A}$ we have%
\begin{equation*}
\mathbf{\varphi }_{X}\circ \mathbf{\mu }_{X}=\mathbf{\mu }_{X}^{\prime
}\circ \mathbf{\varphi }_{\mathbb{M}^{\prime }(X)}\circ \mathbb{M}(\mathbf{%
\varphi }_{X})\text{ and }\mathbf{\varphi }_{X}\circ \mathbf{\nu }_{X}=%
\mathbf{\nu }_{X}^{\prime }\circ \mathbf{\xi }_{X}.
\end{equation*}%
For a fixed endo-functor $\mathbb{J}:\mathfrak{A}\longrightarrow \mathfrak{A}%
,$ we denote by $\mathbb{J}$-$\mathbf{Monad}_{\mathfrak{A}}$ the subcategory
of $\mathbf{JMonad}_{\mathfrak{A}}$ of $\mathbb{J}$-monads on $\mathfrak{A}$
with $\mathbf{\omega }$ the identity natural transformation. In the special
case $\mathbb{J}=\mathbb{I}_{\mathfrak{A}}$ and $\mathbf{\omega }$ is the
identity natural transformation, we drop the prefix and recover the
classical notion of \emph{monads} on $\mathfrak{A}.$
\end{punto}

\begin{remark}
\label{J-factor}As we saw above, a $\mathbb{J}$-monad $(\mathbb{M},\mathbf{%
\mu },\mathbf{\omega },\mathbf{\nu };\mathbb{J})$ is a generalized notion of
a monad. However, it can also be seen as just a monad $(\mathbb{M},\mathbf{%
\mu },\mathbf{\eta })$ whose unit $\mathbf{\eta }:=\mathbb{I}\overset{%
\mathbf{\omega }}{\longrightarrow }\mathbb{J}\overset{\mathbf{\nu }}{%
\longrightarrow }\mathbb{M}$ factorizes through $\mathbb{J}.$ Having this in
mind, a morphism $(\mathbf{\varphi };\mathbf{\xi }):(\mathbb{M},\mathbf{\mu }%
,\mathbf{\omega },\mathbf{\nu };\mathbb{J})\longrightarrow (\mathbb{M}%
^{\prime },\mathbf{\mu }^{\prime },\mathbf{\omega }^{\prime },\mathbf{\nu }%
^{\prime };\mathbb{J}^{\prime })$ in $\mathbf{JMonad}_{\mathfrak{A}}$ is
just a morphism of monads which is compatible with the factorizations of the
units through $\mathbb{J}$ and $\mathbb{J}^{\prime }.$
\end{remark}

\begin{punto}
Let $(\mathbb{M},\mathbf{\mu },\mathbf{\omega },\mathbf{\nu };\mathbb{J})\in
\mathbf{JMonad}_{\mathfrak{A}}.$ An $(\mathbb{M};\mathbb{J})$\emph{-module}
is an object $X\in \mathrm{Obj}(\mathfrak{A})$ with a morphism $\varrho _{X}:%
\mathbb{M}(X)\longrightarrow X$ such that the following diagrams are
commutative%
\begin{equation*}
\begin{tabular}{lll}
$\xymatrix{ \mathbb{MM}(X) \ar[dd]_{\mathbf{\mu}_X}
\ar[rr]^{\mathbb{M}(\varrho_{X})} & & \mathbb{M} (X) \ar[dd]^{\varrho_X} \\
\\ \mathbb{M}(X) \ar[rr]_{\varrho_{X}} && X}$ &  & $\xymatrix{\mathbb{M}(X)
\ar[rr]^{\mathbf{\varrho}_X} & & X \ar[dd]^{{\mathbf{\omega}}_X} \\ \\
\mathbb{J}(X) \ar@{=}[rr] \ar[uu]^{{\mathbf{\nu}}_X}& & \mathbb{J}(X)}$%
\end{tabular}%
\end{equation*}%
The category of $(\mathbb{M};\mathbb{J})$-modules and morphisms those of $%
\mathbb{M}$-acts is denoted by $\mathfrak{A}_{(\mathbb{M};\mathbb{J})}.$ In
case $\mathbb{J}\simeq \mathbb{I}_{\mathfrak{A}}$ and $\mathbf{\omega }$ is
the identity natural transformation, we recover the category of $\mathbb{M}$%
\emph{-modules} of the monad $\mathbb{M}.$
\end{punto}

\begin{punto}
Let $(\mathbb{M},\mathbf{\mu },\mathbf{\omega },\mathbf{\nu };\mathbb{J})\in
\mathbf{JMonad}_{\mathfrak{A}}.$ For every $X\in \mathrm{Obj}(\mathbb{A}),$ $%
\mathbb{M}(X)$ is an $(\mathbb{M};\mathbb{J})$-module through%
\begin{equation*}
\varrho _{\mathbb{M}(X)}:\mathbb{M}\mathbf{(}\mathbb{M}(X))\overset{\mathbf{%
\mu }_{X}}{\longrightarrow }\mathbb{M}(X).
\end{equation*}%
Such modules are called \emph{free }$(\mathbb{M};\mathbb{J})$\emph{-modules }%
and we have the so called \emph{free functor}%
\begin{equation*}
\mathcal{F}_{(\mathbb{M};\mathbb{J})}:\mathfrak{A}\longrightarrow \mathfrak{A%
}_{(\mathbb{M};\mathbb{J})},\text{ }X\mapsto \mathbb{M}(X).
\end{equation*}%
The full subcategory of free $(\mathbb{M};\mathbb{J})$-modules is called the
\emph{Kleisli category} and is denoted by $\widetilde{\mathfrak{A}}_{(%
\mathbb{M};\mathbb{J})}.$
\end{punto}

\begin{remark}
\label{J-lifted}Let $(\mathbb{M},\mathbf{\mu },\mathbf{\omega },\mathbf{\nu }%
;\mathbb{J})\in \mathbf{JMonad}_{\mathfrak{A}}$ with $\mathbb{MJ}\simeq
\mathbb{JM}.$ If $X$ is an $(\mathbb{M};\mathbb{J})$-module, then $\mathbb{J}%
(X)$ is also an $(\mathbb{M};\mathbb{J})$-module through%
\begin{equation*}
\varrho _{\mathbb{J}(X)}:\mathbb{MJ}(X)\simeq \mathbb{JM}(X)\overset{\mathbb{%
J}(\varrho _{X})}{\longrightarrow }\mathbb{J}(X).
\end{equation*}%
Moreover, if $Y=\mathbb{M}(X)$ is a free $(\mathbb{M};\mathbb{J})$-module,
then $\mathbb{J}(Y)=\mathbb{JM}(X)\simeq \mathbb{MJ}(X)$ is also a free $(%
\mathbb{M};\mathbb{J})$-module. One can easily see that $\mathbb{J}$ can be
lifted to endo-functors $\mathbb{J}^{\prime }:$ $\mathfrak{A}_{(\mathbb{M};%
\mathbb{J})}\longrightarrow \mathfrak{A}_{(\mathbb{M};\mathbb{J})}$ and $%
\widetilde{\mathbb{J}}:$ $\widetilde{\mathfrak{A}}_{(\mathbb{M};\mathbb{J}%
)}\longrightarrow \widetilde{\mathfrak{A}}_{(\mathbb{M};\mathbb{J})}.$
\end{remark}

\begin{punto}
Let $(\mathbb{M},\mathbf{\mu },\mathbf{\omega },\mathbf{\nu };\mathbb{J})\in
\mathbf{JMonad}_{\mathfrak{A}}$ and assume that $\mathbb{MJ}\simeq \mathbb{JM%
}.$ Then we have a natural isomorphism for every $X\in \mathfrak{A}$ and $%
Y\in \mathfrak{A}_{(\mathbb{M};\mathbb{J})}:$%
\begin{equation*}
\mathfrak{A}_{(\mathbb{M};\mathbb{J})}(\mathcal{F}_{(\mathbb{M};\mathbb{J}%
)}(X),\mathbb{J}(Y))\simeq \mathfrak{A}(X,\mathbb{J}(Y)),\text{ }f\mapsto
f\circ (\mathbf{\nu }\circ \mathbf{\omega })_{X}
\end{equation*}%
with inverse $g\longmapsto \varrho _{\mathbb{J}(Y)}\circ \mathcal{F}_{(%
\mathbb{M};\mathbb{J})}(g).$ Consider the forgetful functor $U:\mathfrak{A}%
_{(\mathbb{M};\mathbb{J})}\longrightarrow \mathfrak{A}$ and the endo-functor
$\mathbb{J}^{\prime }:$ $\mathfrak{A}_{(\mathbb{M};\mathbb{J}%
)}\longrightarrow \mathfrak{A}_{(\mathbb{M};\mathbb{J})}$ (see Remark \ref%
{J-lifted}). Then we have a natural isomorphism%
\begin{equation*}
\mathfrak{A}_{(\mathbb{M};\mathbb{J})}(\mathcal{F}_{(\mathbb{M};\mathbb{J})}(%
\mathbb{J}(X)),\mathbb{J}^{\prime }(Y))\simeq \mathfrak{A}(\mathbb{J}(X),U(%
\mathbb{J}^{\prime }(Y)));
\end{equation*}%
\emph{i.e.} $(\mathcal{F}_{(\mathbb{M};\mathbb{J})}(-),U)$ is a $(\mathbb{J},%
\mathbb{J}^{\prime })$-adjoint pair.
\end{punto}

\subsection*{$\mathbb{J}$-Comonads}

\begin{punto}
Let $\mathbb{J}$ be an endo-functor on $\mathfrak{A}.$ With a $\mathbb{J}$-%
\emph{comonad} on $\mathfrak{A}$ we mean a datum $(\mathbb{C},\mathbf{\Delta
},\mathbf{\omega },\mathbf{\theta })$ consisting of an endo-functor $\mathbb{%
C}:\mathfrak{A}\longrightarrow \mathfrak{A}$ associated with natural
transformations%
\begin{equation*}
\mathbf{\Delta }:\mathbb{C}\longrightarrow \mathbb{CC},\text{ }\mathbf{%
\omega }:\mathbb{I}\longrightarrow \mathbb{J}\text{ and }\mathbf{\theta }:%
\mathbb{C}\longrightarrow \mathbb{J}
\end{equation*}%
such that the following diagrams are commutative
\begin{equation*}
\begin{tabular}{lll}
$\xymatrix{\mathbb{C} \ar[dd]_{\mathbf{\Delta}} \ar[rr]^{\mathbf{\Delta}} &
& \mathbb{CC} \ar[dd]^{\mathbf{\Delta}\mathbb{C}} \\ & & \\ \mathbb{CC}
\ar[rr]_{\mathbb{C}\mathbf{\Delta} } & & \mathbb{CCC}}$ & $%
\xymatrix{\mathbb{I C} \ar[dd]_{\mathbf{\omega}\mathbb{C}}
\ar[rr]^{\mathbf{\Delta}} & & \mathbb{CC}
\ar[dd]^{\mathbf{\theta}\mathbb{C}} \\ & & \\ \mathbb{J C} \ar@{=}[rr] & &
\mathbb{J C}}$ & $\xymatrix{\mathbb{C I} \ar[dd]_{\mathbb{C}\mathbf{\omega}}
\ar[rr]^{\mathbf{\Delta}} & & \mathbb{CC}
\ar[dd]^{\mathbb{C}\mathbf{\theta}} \\ & & \\ \mathbb{C J} \ar@{=}[rr] & &
\mathbb{C J}}$%
\end{tabular}%
\end{equation*}%
\emph{i.e. }for every $X\in \mathfrak{A}$ we have%
\begin{equation*}
\mathbf{\Delta }_{\mathbb{C}(X)}\circ \mathbf{\Delta }_{X}=\mathbb{C}(%
\mathbf{\Delta }_{X})\circ \mathbf{\Delta }_{X},\text{ }\mathbf{\theta }_{%
\mathbb{C}(X)}\circ \mathbf{\Delta }_{X}=\mathbf{\omega }_{\mathbb{C}(X)}%
\text{ and }\mathbb{C}(\mathbf{\theta }_{X})\circ \mathbf{\Delta }_{X}=%
\mathbb{C}(\mathbf{\omega }_{X}).
\end{equation*}
\end{punto}

\begin{punto}
By $\mathbf{JComonad}_{\mathfrak{A}}$ we denote the category whose objects
are $\mathbb{J}$-comonads, where $\mathbb{J}$ runs over the class of
endo-functors on $\mathfrak{A}.$ A morphism $(\mathbf{\psi };\mathbf{\xi }):(%
\mathbb{C},\mathbf{\Delta },\mathbf{\omega },\mathbf{\theta };\mathbb{J}%
)\longrightarrow (\mathbb{C}^{\prime },\mathbf{\Delta }^{\prime },\mathbf{%
\omega }^{\prime },\mathbf{\theta }^{\prime };\mathbb{J}^{\prime })$ in this
category consists of natural transformations $\mathbf{\psi }:\mathbb{C}%
\longrightarrow \mathbb{C}^{\prime }$ and $\mathbf{\xi }:\mathbb{J}%
\longrightarrow \mathbb{J}^{\prime }$ such that the following diagrams are
commutative%
\begin{equation*}
\begin{tabular}{lll}
$\xymatrix{ \mathbb{C} \ar[d]_{\mathbf{\psi}} \ar[rr]^{\mathbf{\Delta}} & &
\mathbb{CC} \ar[d]^{\mathbf{ \psi \psi}} \\ {\mathbb{C'}}
\ar[rr]_{{\mathbf{\Delta}}'} & & \mathbb{C'C'} }$ &  & $\xymatrix{\mathbb{C}
\ar[rr]^{\mathbf{\theta}} \ar[d]_{\mathbf{\psi}} & & \mathbb{J}
\ar[d]^{\mathbf{\xi}} \\ \mathbb{C}' \ar[rr]_{\mathbf{\theta}'} & &
\mathbb{J}'}$%
\end{tabular}%
\end{equation*}%
\emph{i.e.} for every $X\in \mathfrak{A}$ we have%
\begin{equation*}
\mathbf{\psi }_{\mathbb{C}^{\prime }(X)}\circ \mathbb{C}(\mathbf{\psi }%
_{X})\circ \mathbf{\Delta }_{X}=\mathbf{\Delta }_{X}^{\prime }\circ \mathbf{%
\psi }_{X}\text{ and }\mathbf{\xi }_{X}\circ \mathbf{\theta }_{X}=\mathbf{%
\theta }_{X}^{\prime }\circ \mathbf{\psi }_{X}.
\end{equation*}%
For a fixed endo-functor $\mathbb{J}:\mathfrak{A}\longrightarrow \mathfrak{A}%
,$ we denote by $\mathbb{J}$-$\mathbf{Comonad}_{\mathfrak{A}}$ the
subcategory of $\mathbb{J}$-comonads on $\mathfrak{A}$ with $\mathbf{\omega }
$ the identity transformation. When we drop the prefix, we have the special
case $\mathbb{J}=\mathbb{I}_{\mathfrak{A}}$ and $\mathbf{\omega }$ is the
identity natural transformation and recover the notion of \emph{comonads} on
$\mathfrak{A}.$
\end{punto}

\begin{remark}
$\mathbb{J}$-Comonads are \emph{not} fully dual to $\mathbb{J}$-monads.
Recall from Remark \ref{J-factor} that a $\mathbb{J}$-monad can be seen as a
monad whose unit factorizes through $\mathbb{J}.$ On the other hand, $%
\mathbb{J}$-comonads cannot be seen as a special type of comonads. The lack
of duality is because not all arrows are reversed; the arrow $\mathbf{\omega
}:\mathbb{I}\longrightarrow \mathbb{J}$ is assumed for both. Notice that
keeping this arrow is suggested by the concrete example in Section 5.
\end{remark}

\begin{punto}
Let $(\mathbb{C},\mathbf{\Delta },\mathbf{\omega ,\theta };\mathbb{J})\in
\mathbf{JComonad}_{\mathfrak{A}}.$ A $(\mathbb{C};\mathbb{J})$\emph{-comodule%
} is an object $X\in \mathrm{Obj}(\mathfrak{A})$ along with a morphism $%
\varrho ^{X}:X\longrightarrow \mathbb{C}(X)$ in $\mathfrak{A}$ such that the
following diagrams are commutative%
\begin{equation*}
\begin{tabular}{lll}
$\xymatrix{X \ar[dd]_{\varrho^X} \ar[rr]^{\varrho^{X}} & & \mathbb{C}(X)
\ar[dd]^{\mathbb{C}(\varrho^X)} \\ & & \\ \mathbb{C}(X)
\ar[rr]_{\mathbf{\Delta}_X} & & \mathbb{CC}(X)}$ &  & $\xymatrix{X
\ar[dd]_{{\mathbf{\omega}_X}} \ar[rr]^{\varrho ^X} & & \mathbb{C}(X)
\ar[dd]^{\mathbf{\theta}_X}\\ & & \\ \mathbb{J}(X) \ar@{=}[rr] & &
\mathbb{J}(X)}$%
\end{tabular}%
\end{equation*}%
The category of $(\mathbb{C};\mathbb{J})$-comodules and morphisms those of $%
\mathbb{C}$-coacts is denoted by $\mathfrak{A}^{(\mathbb{C};\mathbb{J})}.$
In case $\mathbb{J}=\mathbb{I}_{\mathfrak{A}}$ and $\mathbf{\omega }$ is the
identity natural transformation, we recover the category of $\mathbb{C}$%
\emph{-comodules} for the comonad $\mathbb{C}.$
\end{punto}

\begin{punto}
Let $(\mathbb{C},\mathbf{\Delta },\mathbf{\varepsilon };\mathbb{J})\in
\mathbf{JComonad}_{\mathfrak{A}}.$ For every $X\in \mathrm{Obj}(\mathfrak{A}%
),$ $\mathbb{C}(X)$ has a canonical structure of a $(\mathbb{C};\mathbb{J})$%
-comodule through%
\begin{equation*}
\varrho ^{\mathbb{C}(X)}:\mathbb{C}(X)\overset{\mathbf{\Delta }_{X}}{%
\longrightarrow }\mathbb{CC}(X).
\end{equation*}%
Such comodules are called \emph{cofree }$(\mathbb{C};\mathbb{J})$\emph{%
-comodules }and we have the so called \emph{cofree functor}%
\begin{equation*}
\mathcal{F}^{\mathbb{C}}:\mathfrak{A}\longrightarrow \mathfrak{A}^{(\mathbb{C%
};\mathbb{J})},\text{ }X\mapsto \mathbb{C}(X).
\end{equation*}%
The full subcategory of cofree $(\mathbb{C};\mathbb{J})$-comodules is called
the \emph{Kleisli category} of $\mathbb{C}$ and is denoted by $\widetilde{%
\mathfrak{A}}^{(\mathbb{C};\mathbb{J})}.$
\end{punto}

\begin{remark}
Let $(\mathbb{C},\mathbf{\Delta },\mathbf{\omega },\mathbf{\theta };\mathbb{J%
})\in \mathbf{JComonad}_{\mathfrak{A}}$ with $\mathbb{JC}\simeq \mathbb{CJ}.$
If $X$ is a $(\mathbb{C};\mathbb{J})$-comodule, then $\mathbb{J}(X)$ is also
a $(\mathbb{C};\mathbb{J})$-comodule through%
\begin{equation*}
\varrho ^{\mathbb{J}(X)}:\mathbb{J}(X)\overset{\mathbb{J}(\varrho ^{X})}{%
\longrightarrow }\mathbb{JC}(X)\simeq \mathbb{C}(\mathbb{J}(X)).
\end{equation*}%
If $Y=\mathbb{C}(X)$ is a cofree $(\mathbb{C};\mathbb{J})$-comodule, then $%
\mathbb{J}(Y)=\mathbb{JC}(X)\simeq \mathbb{CJ}(X)$ is also a cofree $(%
\mathbb{C},\mathbb{J})$-comodule. One case easily see that $\mathbb{J}$
lifts to endo-functors $\mathbb{J}^{\prime }:\mathfrak{A}^{(\mathbb{C};%
\mathbb{J})}\longrightarrow \mathfrak{A}^{(\mathbb{C};\mathbb{J})}$ and $%
\widetilde{\mathbb{J}}:\widetilde{\mathfrak{A}}^{(\mathbb{C};\mathbb{J}%
)}\longrightarrow \widetilde{\mathfrak{A}}^{(\mathbb{C};\mathbb{J})}.$
\end{remark}

\begin{punto}
Let $(\mathbb{C},\mathbf{\Delta },\mathbf{\omega },\mathbf{\theta };\mathbb{J%
})\in \mathbf{JComonad}_{\mathfrak{A}}$ with $\mathbb{J}$ idempotent and $%
\mathbb{JC}\simeq \mathbb{CJ}.$ Consider the forgetful functor $U:\mathfrak{A%
}^{(\mathbb{C};\mathbb{J})}\longrightarrow \mathfrak{A}$ and the
endo-functor $\mathbb{J}^{\prime }:\mathfrak{A}^{(\mathbb{C};\mathbb{J}%
)}\longrightarrow \mathfrak{A}^{(\mathbb{C};\mathbb{J})}.$ Then we have a
natural isomorphism for $X\in \mathfrak{A}$ and $Y\in \mathfrak{A}^{(\mathbb{%
C};\mathbb{J})}:$%
\begin{equation*}
\mathfrak{A}^{(\mathbb{C};\mathbb{J})}(\mathbb{J}^{\prime }(Y),\mathcal{F}^{%
\mathbb{C}}(\mathbb{J}(X)))\simeq \mathfrak{A}(U(\mathbb{J}^{\prime }(Y)),%
\mathbb{J}(X)),\text{ }f\mapsto \mathbf{\theta }_{\mathbb{J}(X)}\circ f
\end{equation*}%
with inverse $g\longmapsto \mathcal{F}^{\mathbb{C}}(g)\circ \varrho ^{%
\mathbb{J}(Y)};$ \emph{i.e.} $(U,\mathcal{F}^{(\mathbb{C};\mathbb{J})}(-))$
is a $(\mathbb{J}^{\prime },\mathbb{J})$-adjoint pair.
\end{punto}

\begin{proposition}
Let $\mathfrak{A}$ and $\mathfrak{B}$ be categories, $L:\mathfrak{A}%
\longrightarrow \mathfrak{B},$ $R:\mathfrak{B}\longrightarrow \mathfrak{A}$
be functors and $\mathbb{J}:\mathfrak{A}\longrightarrow \mathfrak{A},$ $%
\mathbb{K}:\mathfrak{B}\longrightarrow \mathfrak{B}$ endo-functors such that
$L\mathbb{J}\simeq \mathbb{K}L,$ $\mathbb{J}R\simeq R\mathbb{K}$ and $(L,R)$
is a $(\mathbb{J},\mathbb{K})$-adjoint pair.

\begin{enumerate}
\item $(\mathcal{L},R)$ is an adjoint pair where $\mathcal{L}:\mathbb{J}(%
\mathfrak{A})\overset{L}{\longrightarrow }\mathbb{K}(\mathfrak{B})$ and $%
\mathcal{R}:\mathbb{K}(\mathfrak{B})\overset{R}{\longrightarrow }\mathbb{J}(%
\mathfrak{A})$ with unit and counit of adjunction given by%
\begin{equation*}
\eta :\mathbb{J}\longrightarrow RL\mathbb{J}\text{ and }\varepsilon :LR%
\mathbb{K}\longrightarrow \mathbb{K}.
\end{equation*}

\item $RL$ is a monad on $\mathbb{J}(\mathfrak{A})$ with%
\begin{equation*}
\mathbf{\mu }_{RL}:(RL\mathbf{)(}RL\mathbf{)}\mathbb{J}\mathbf{\simeq }R%
\mathbf{(}LR\mathbb{K}\mathbf{)}L\overset{R\mathbf{\varepsilon }L}{\mathbf{%
\longrightarrow }}R\mathbb{K}L\simeq (RL\mathbf{)}\mathbb{J}\text{ and }%
\mathbf{\eta }_{RL}:=\eta .
\end{equation*}

\item $LR$ is a comonad on $\mathbb{K}(\mathfrak{B})$ with%
\begin{equation*}
\Delta _{LR}:(LR)\mathbb{K}\simeq L\mathbb{J}R\overset{L\mathbf{\eta }R}{%
\longrightarrow }L(RL\mathbb{J})R\simeq (LR)(LR)\mathbb{K}\text{ and }%
\mathbf{\varepsilon }_{LR}:=\varepsilon .
\end{equation*}

\item $L$ is a monad on $\mathbb{J}(\mathfrak{A})$ if and only if $R$ is a
comonad on $\mathbb{K}(\mathfrak{B}).$ In this case, $\mathbb{J}(\mathfrak{A}%
)_{L}\simeq \mathbb{K}(\mathfrak{B})^{R}.$

\item $L$ is a comonad on $\mathbb{J}(\mathfrak{A})$ if and only if $R$ is a
monad on $\mathbb{K}(\mathfrak{B}).$ In this case, $\widetilde{\mathbb{J}(%
\mathfrak{A})}^{L}\simeq \widetilde{\mathbb{K}(\mathfrak{B})}_{R}.$
\end{enumerate}
\end{proposition}

\begin{Beweis}
By assumption $L\mathbb{J}\simeq \mathbb{K}L$ whence $\mathcal{L}(\mathbb{J}(%
\mathfrak{A})):=L\mathbb{J}(\mathfrak{A})=\mathbb{K}L(\mathfrak{A})\subseteq
\mathbb{K}(\mathfrak{B})$ and $\mathbb{J}R\simeq R\mathbb{K}$ whence $%
\mathcal{R}(\mathbb{K}(\mathfrak{B})):=R(\mathbb{K}(\mathfrak{B}))=\mathbb{J}%
R(\mathfrak{B})\subseteq \mathbb{J}(\mathfrak{A}).$ The assumptions imply
that $(\mathcal{L},\mathcal{R})$ is an adjoint pair. The result follows now
from the classical result on right adjoint pairs (\emph{e.g.} \cite[%
Proposition 3.1]{EM1965}, \cite[2.5, 2.6]{BBW2009})\emph{.}$\blacksquare $
\end{Beweis}

\section{Semiunital Semimonoidal Categories}

\qquad A semimonoidal category is roughly speaking a monoidal category not
necessarily with a unit object. The reader might consult the literature for
the precise definitions and for the notions of (op)-semimonoidal functors
between such categories. In this section, we introduce a notion of \emph{%
semiunital semimonoidal categories} and \emph{semiunital} (\emph{op}-)\emph{%
semimonoidal functors}.

\subsection*{Semiunits}

\begin{punto}
Let $(\mathcal{V},\bullet )$ be a semimonoidal category. We say that $%
\mathbf{I}\in \mathcal{V}$ is a \emph{semiunit} iff

\begin{enumerate}
\item there is a natural transformation $\mathbf{\omega }:\mathbb{I}%
\longrightarrow (\mathbf{I}\bullet -);$

\item there exists an isomorphisms of functors $\mathbf{I}\bullet -\simeq
-\bullet \mathbf{I},$ \emph{i.e.} there is a natural isomorphism $\mathbf{I}%
\bullet X\overset{\ell _{X}}{\simeq }X\bullet \mathbf{I}$ in $\mathcal{V}$
with inverse $\wp _{X},$ for each object $X$ of $\mathcal{V},$ such that $%
\ell _{\mathbf{I}}=\wp _{\mathbf{I}}$ and the following diagrams are
commutative for all $X,Y\in \mathcal{V}:$%
\begin{equation*}
\begin{tabular}{l}
$\xymatrix{(\mathbf{I} \bullet X) \bullet Y
\ar[rr]^{{\mathbf{\gamma}}_{\mathbf{I},X,Y}} \ar[dd]_{\ell_{X} \bullet Y} &
& \mathbf{I} \bullet (X \bullet Y) \ar[rr]^{\ell_{X \bullet Y}} & & (X
\bullet Y) \bullet \mathbf{I} \ar[dd]^{{\mathbf{\gamma}}_{X,Y,\mathbf{I}}}
\\ & & & & & \\ (X \bullet \mathbf{I}) \bullet Y
\ar[rr]_{{\mathbf{\gamma}}_{X,\mathbf{I},Y}} & & X \bullet (\mathbf{I}
\bullet Y) \ar[rr]_{X \bullet \ell_{Y}} & & X \bullet (Y \bullet \mathbf{I})
}$ \\
\\
$\xymatrix{(\mathbf{I} \bullet X) \bullet Y \ar[rrdd]_{\simeq} & & X \bullet
Y \ar[ll]_{{{{\mathbf{\omega}}_{X}} \bullet Y}}
\ar[dd]_{{{\mathbf{\omega}}_{X \bullet Y}}} \ar[rr]^{X \bullet
{{{\mathbf{\omega}}_{Y}}}} & & X \bullet (\mathbf{I} \bullet Y)
\ar[lldd]^{\simeq} \\ & & & & & \\ & & \mathbf{I} \bullet (X \bullet Y) & & }
$%
\end{tabular}%
\end{equation*}%
If $X\overset{\mathbf{\omega }_{X}}{\simeq }\mathbf{I}\bullet X$ ($\overset{%
\ell _{X}}{\simeq }X\bullet \mathbf{I}$), then we say that $X$ is \emph{firm}%
.
\end{enumerate}
\end{punto}

\begin{notation}
If $X$ is firm, then we set $\lambda _{X}:=\mathbf{\omega }_{X}^{-1}:\mathbf{%
I}\bullet X\longrightarrow X$ and $\rho _{X}:X\bullet \mathbf{I}\overset{\wp
_{X}}{\simeq }\mathbf{I}\bullet X\overset{\mathbf{\omega }_{X}^{-1}}{%
\longrightarrow }X.$ With $\mathcal{V}^{\mathrm{firm}}$ we denote the \emph{%
full} subcategory of firm objects in $\mathcal{V}.$
\end{notation}

\begin{remark}
If $\mathbf{I}$ is firm (called also \emph{pseudo-idempotent}) and $\mathbf{%
\omega }_{\mathbf{I}}^{-1}\bullet \mathbf{I}=\mathbf{I\bullet \omega
_{I}^{-1}},$ then one says that $\mathbf{I}$ is \emph{idempotent} \cite%
{Koc2008}.
\end{remark}

\begin{remark}
Let $(\mathcal{V},\bullet )$ be a semimonoidal category. One says that $%
\mathcal{V}$ is \emph{monoidal} \cite{Mac1998} iff $\mathcal{V}$ has a \emph{%
unit} (or an $LR$ \emph{unit}), \emph{i.e.} a distinguished object $\mathbf{I%
}\in \mathcal{V}$ with natural isomorphisms $\mathbf{I}\bullet X\overset{%
\lambda _{X}}{\simeq }X$ and $X\bullet \mathbf{I}\overset{\rho _{X}}{\simeq }%
X$ such that $X\bullet \lambda _{Y}=\rho _{X}\bullet Y$ for all $X,Y\in
\mathcal{V}$ (equivalently, $\lambda _{\mathbf{I}}=\rho _{\mathbf{I}},$ $%
\lambda _{X\bullet Y}=\lambda _{X}\bullet Y$ and $\rho _{X\bullet
Y}=X\bullet \rho _{Y}$ for all $X,Y\in \mathcal{V}$). Kock \cite{Koc2008}
called an object $\mathbf{I}\in \mathcal{V}$ a \emph{Saavedra unit} --
called also a \emph{reduced unit} -- iff it is pseudo-idempotent and \emph{%
cancellable} in the sense that the endo-functors $\mathbf{I}\bullet -$ and $%
-\bullet \mathbf{I}$ are full and faithful (equivalently, $\mathbf{I}$ is
idempotent and the endo-functors $\mathbf{I}\bullet -$ and $-\bullet \mathbf{%
I}$ are equivalences of categories). Moreover, he showed that $\mathbf{I}$
is a unit if and only if $\mathbf{I}$ is a Saavedra unit. Indeed, every unit
is a semiunit, whence our notion of semiunital semimonoidal categories
generalizes the classical notion of monoidal categories.
\end{remark}

\begin{punto}
Let $(\mathcal{V},\bullet ,\mathbf{I}_{\mathcal{V}};\mathbf{\omega }_{%
\mathcal{V}})$ and $(\mathcal{W},\otimes ,\mathbf{I}_{\mathcal{W}};\mathbf{%
\omega }_{\mathcal{W}})$ be semiunital semimonoidal categories. A
semimonoidal functor $F:\mathcal{V}\longrightarrow \mathcal{W},$ with a
natural transformation ${\mathbf{\phi }}:F(-)\otimes F(-)\longrightarrow
F(-\bullet -),$ is said to be \emph{semiunital semimonoidal }iff there
exists a \emph{coherence morphism} $\phi :\mathbf{I}_{\mathcal{W}%
}\longrightarrow F(\mathbf{I}_{\mathcal{V}})$ in $\mathcal{W}$ such that the
following diagram is commutative%
\begin{equation*}
\begin{tabular}{l}
$\xymatrix{\mathbf{I}_\W \otimes F(X) \ar[rr]^{\ell_{F(X)}} \ar[dd]_{\phi
\otimes F(X)} & & F(X)\otimes \mathbf{I}_\W \ar[dd]^{F(X) \otimes \phi} \\ &
& \\ F(\mathbf{I}_{\V}) \otimes F(X) \ar[dd]_{{\mathbf{\phi}}_{\mathbf{I}_\V
,X}} & & F(X)\otimes
F(\mathbf{I}_\V)\ar[dd]^{{\mathbf{\phi}}_{X,\mathbf{I}_\V}} \\ & & \\
F(\mathbf{I}_\V \bullet X) \ar[rr]_{F(\ell_X)} & & F(X\bullet \mathbf{I}_\V)
}$%
\end{tabular}%
\end{equation*}%
Moreover, we say that $F$ is a \emph{strong} (\emph{strict}) \emph{%
semiunital semimonoidal functor} iff $F$ is strong (strict) as a
semimonoidal functor and ${\mathbf{\phi }}$ is an isomorphism (identity).
For two semimonoidal functors $F,F^{\prime }:\mathcal{V}\longrightarrow
\mathcal{W},$ we say that a semimonoidal natural transformation $\mathbf{\xi
}:F\longrightarrow F^{\prime }$ is \emph{semiunital semimonoidal} iff the
following diagram is commutative%
\begin{equation*}
\xymatrix{ & \mathbf{I}_{\mathcal{W}} \ar[dl]_{\phi} \ar[dr]^{\phi'} & \\
F(\mathbf{I}_{\mathcal{V}})
\ar[rr]_{{\mathbf{\xi}}_{\mathbf{I}_{\mathcal{V}}}} & &
F'(\mathbf{I}_{\mathcal{V}}) }
\end{equation*}%
One can dually define \emph{semiunital }(\emph{strong}, \emph{strict}) \emph{%
op-semimonoidal functors }and semiunital natural transformations between
them.
\end{punto}

\begin{remarks}
Let $(\mathcal{V},\bullet ,\mathbf{I};\mathbf{\omega })$ be a semiunital
semimonoidal category and consider the functor%
\begin{equation*}
\mathbb{J}:=\mathbf{I}\bullet -:\mathcal{V}\longrightarrow \mathcal{V}.
\end{equation*}

\begin{enumerate}
\item We have natural isomorphisms%
\begin{equation*}
\mathbb{J}(\mathbf{I})\bullet X\simeq \mathbb{J}(\mathbb{J}(X))\simeq
X\bullet \mathbb{J}(\mathbf{I})\text{ and }\mathbb{J}(X)\bullet Y\simeq
X\bullet \mathbb{J}(Y)
\end{equation*}%
for any $X,Y\in \mathcal{V}.$

\item $\mathbb{J}$ is op-semimonoidal; the natural transformation $\mathbf{%
\delta }_{X,Y}:\mathbb{J}(-\bullet -)\longrightarrow \mathbb{J}(-)\bullet
\mathbb{J}(-)$ is given by the collection of morphisms%
\begin{equation*}
{\mathbf{\delta }}_{X,Y}:\mathbf{I}\bullet (X\bullet Y)\overset{\mathbf{%
\omega }_{\mathbf{I}}\bullet (X\bullet Y)}{\longrightarrow }(\mathbf{I}%
\bullet \mathbf{I})\bullet (X\bullet Y)\simeq (\mathbf{I}\bullet X)\bullet (%
\mathbf{I}\bullet Y)
\end{equation*}%
for all $X,Y\in \mathcal{V}.$

\item Assume that $\mathbf{I}$ is firm.

\begin{enumerate}
\item $\mathbb{J}$ is strong semiunital semimonoidal with%
\begin{equation*}
{\mathbf{\phi }}_{X,Y}:(\mathbf{I}\bullet X)\bullet (\mathbf{I}\bullet
Y)\simeq (\mathbf{I}\bullet \mathbf{I})\bullet (X\bullet Y)\overset{\mathbf{%
\omega }_{\mathbf{I}}^{-1}\bullet (X\bullet Y)}{\simeq }\mathbf{I}\bullet
(X\bullet Y)\text{ and }\phi :=\omega _{\mathbf{I}}:\mathbf{I}%
\longrightarrow \mathbf{I}\bullet \mathbf{I}
\end{equation*}%
for all $X,Y\in \mathcal{V}.$

\item $\mathbb{J}$ is strong semiunital op-semimonoidal with%
\begin{equation*}
\delta :=\omega _{\mathbf{I}}^{-1}:\mathbf{I}\bullet \mathbf{I}%
\longrightarrow \mathbf{I}.
\end{equation*}

\item the full subcategory $(\mathbf{U}\mathcal{V},\bullet ,\mathbf{I})$ is
monoidal.

\item $(\mathbb{J}(\mathcal{V}),\bullet ,\mathbf{I})$ is a monoidal full
subcategory of $(\mathbf{U}\mathcal{V},\bullet ,\mathbf{I})$ with%
\begin{eqnarray*}
\lambda _{\mathbf{I}\bullet X} &:&\mathbf{I}\bullet (\mathbf{I}\bullet X)%
\overset{{\mathbf{\gamma }}_{\mathbf{I},\mathbf{I},X}^{-1}}{\simeq }(\mathbf{%
I}\bullet \mathbf{I})\bullet X\overset{\mathbf{\omega }_{\mathbf{I}}^{-1}}{%
\simeq }\mathbf{I}\bullet X; \\
\rho _{\mathbf{I}\bullet X} &:&(\mathbf{I}\bullet X)\bullet \mathbf{I}%
\overset{{\mathbf{\gamma }}_{\mathbf{I},X,\mathbf{I}}}{\simeq }\mathbf{I}%
\bullet (X\bullet \mathbf{I})\overset{\mathbf{I}\bullet \wp _{X}}{\simeq }%
\mathbf{I}\bullet (\mathbf{I}\bullet X)\overset{\lambda _{\mathbf{I}\bullet
X}}{\simeq }\mathbf{I}\bullet X
\end{eqnarray*}%
for every $X\in \mathcal{V}.$
\end{enumerate}
\end{enumerate}
\end{remarks}

\begin{definition}
Let $(\mathcal{V},\bullet ,\mathbf{I};\mathbf{\omega })$ be a semiunital
semimonoidal category. We say that $V\in \mathcal{V}$ has a \emph{left dual}
iff there exists $V_{l}^{\circledast }\in \mathcal{V}$ along with morphisms $%
\upsilon :\mathbf{I}\longrightarrow \mathbf{I}\bullet V\bullet
V_{l}^{\circledast }$ and $\varpi :\mathbf{I}\bullet V_{l}^{\circledast
}\bullet V\longrightarrow \mathbf{I}$ in $\mathcal{V}$ such that%
\begin{equation*}
(V\bullet \varpi )\circ (\ell _{V}\bullet V_{l}^{\circledast }\bullet
V)\circ (\upsilon \bullet V)=\ell _{V}\text{ and }(\varpi \bullet
V_{l}^{\circledast })\circ (\wp _{V^{\circledast }}\bullet V\bullet
V_{l}^{\circledast })\circ (V_{l}^{\circledast }\bullet \upsilon )=\wp
_{V^{\circledast }}.
\end{equation*}%
A \emph{right dual }$V_{r}^{\circledast }$ of $V$ is defined symmetrically.
We say that $\mathcal{V}$ is \emph{left} (\emph{right}) \emph{autonomous},
or \emph{left} (\emph{right}) \emph{rigid} iff every object in $\mathcal{V}$
has a left (right) dual.
\end{definition}

\begin{definition}
Let $(\mathcal{V},\bullet ,\mathbf{I};\mathbf{\omega })$ be a semiunital
semimonoidal category. We say that $\mathcal{V}$ is \emph{right (left) closed%
} iff for every $V\in \mathcal{V},$ the functor $-\bullet V:\mathbb{J}(%
\mathcal{V})\longrightarrow \mathbb{J}(\mathcal{V})$ ($V\bullet -:\mathbb{J}(%
\mathcal{V})\longrightarrow \mathbb{J}(\mathcal{V})$) has a right-adjoint,
\emph{i.e.} there exists a functor $G:\mathbb{J}(\mathcal{V})\longrightarrow
\mathbb{J}(\mathcal{V})$ and a natural isomorphism for every pair of objects
$X,Y\in \mathcal{V}:$%
\begin{equation*}
\mathcal{V}(X\bullet \mathbf{I}\bullet V,Y\bullet \mathbf{I})\simeq \mathcal{%
V}(X\bullet \mathbf{I},G(Y\bullet \mathbf{I}))\text{ (resp. }\mathcal{V}%
(V\bullet \mathbf{I}\bullet X,Y\bullet \mathbf{I})\simeq \mathcal{V}%
(X\bullet \mathbf{I},G(Y\bullet \mathbf{I}))\text{).}
\end{equation*}%
Moreover, $\mathcal{V}$ is said to be \emph{closed} iff $\mathcal{V}$ is
left and right closed.
\end{definition}

\begin{lemma}
\label{aut->closed}Let $(\mathcal{V},\bullet ,\mathbf{I};\mathbf{\omega })$
be a semiunital semimonoidal category. If $V\in \mathcal{V}$ has a left
\emph{(}right\emph{)} dual $V^{\circledast },$ then $(-\bullet V,-\bullet
V^{\circledast })$ \emph{(}$(V\bullet -,V^{\circledast }\bullet -)$\emph{)}
is a $(\mathbb{J},\mathbb{J})$-adjoint pair. In particular, if $\mathcal{V}$
is left \emph{(}right\emph{)} autonomous, then $\mathcal{V}$ is right \emph{(%
}left\emph{)} closed.
\end{lemma}

\begin{Beweis}
Assume that $V\in \mathcal{V}$ has a left dual $V^{\circledast }.$ For any $%
X,Y\in \mathcal{V}$ we have a natural isomorphism%
\begin{equation}
\mathcal{V}(X\bullet \mathbf{I}\bullet V,Y\bullet \mathbf{I})\simeq \mathcal{%
V}(X\bullet \mathbf{I},Y\bullet \mathbf{I}\bullet V^{\circledast }),\text{ }%
f\longmapsto (f\bullet V^{\circledast })\circ (X\bullet \upsilon )
\label{l-dual-adj}
\end{equation}%
with inverse $g\longmapsto (Y\bullet \varpi )\circ (g\bullet V).\blacksquare
$
\end{Beweis}

\section{Semimonoids and Semicomonoids}

\qquad In this section, we introduce notions of semimonoids and
semicomonoids in semiunital semimonoidal categories.

Throughout, $(\mathcal{V},\bullet ,\mathbf{I};\mathbf{\mathbf{\omega }})$ is
a semiunital semimonoidal category, where $\mathbf{I}$ is a semiunit and $%
\mathbf{\mathbf{\omega }}:\mathbb{I}\longrightarrow \mathbb{J}$ is a natural
transformation between the identity functor and $\mathbb{J}:=\mathbf{I}%
\bullet -\simeq -\bullet \mathbf{I}:\mathcal{V}\longrightarrow \mathcal{V}$
(we assume the existence of natural isomorphisms $\mathbf{I}\bullet X\overset%
{\ell _{X}}{\simeq }X\bullet \mathbf{I}$ with inverse $X\bullet \mathbf{I}%
\overset{\wp _{X}:=\ell _{X}^{-1}}{\longrightarrow }\mathbf{I}\bullet X$ for
every $X\in \mathcal{V}$).

\subsection*{Semimonoids}

\begin{punto}
A $\mathcal{V}$\emph{-semimonoid} consists of a datum $(A,\mu ,\eta ),$
where $A\in \mathcal{V}$ and $\mu :A\bullet A\longrightarrow A,$ $\eta :%
\mathbf{I}\longrightarrow A$ are morphisms in $\mathcal{V}$ such that the
following diagrams are commutative%
\begin{equation*}
\begin{array}{ccc}
\xymatrix{A \bullet A \bullet A \ar[rr]^{\mu \bullet A} \ar[dd]_{A \bullet
\mu } & & A \bullet A \ar[dd]^{\mu}\\ & & \\ A \bullet A \ar[rr]_{\mu} & & A
} &  & \xymatrix{A \bullet A \ar[rr]^{\mu} & & A
\ar[dd]^{{\mathbf{\omega}}_A} & & A \bullet A \ar[ll]_{\mu} \\ & & & & & \\
\mathbf{I} \bullet A \ar[uu]^{\eta \bullet A} \ar@{=}[rr] & & \mathbf{I}
\bullet A \ar[rr]_{\ell_A} & & A \bullet \mathbf{I} \ar[uu]_{A \bullet \eta}}%
\end{array}%
\end{equation*}%
If $A\overset{\mathbf{\mathbf{\omega }}_{A}}{\simeq }\mathbf{I}\bullet A,$
then we say that $A$ is a \emph{unital }$\mathcal{V}$\emph{-semimonoid}. A
\emph{morphism of\ }$\mathcal{V}$\emph{-semimonoids} $f:(A,\mu ,\eta
)\longrightarrow (A^{\prime },\mu ^{\prime },\eta ^{\prime })$ is a morphism
in $\mathcal{V}$ such that the following diagrams are commutative%
\begin{equation*}
\begin{array}{ccc}
\xymatrix{A \bullet A \ar[rr]^{\mu} \ar[dd]_{f \bullet f} & & A
\ar[dd]^{f}\\ & & \\ A' \bullet A' \ar[rr]_{{\mu}'} & & A'} &  & %
\xymatrix{\mathbf{I} \ar@{=}[dd] \ar[rr]^{\eta} & & A \ar[dd]^{f} & & \\ \\
\mathbf{I} \ar[rr]_{_{{\eta}'}} & & A' & &}%
\end{array}%
\end{equation*}%
The category of $\mathcal{V}$-semimonoids is denoted by $\mathbf{SMonoid}(%
\mathcal{V});$ the full subcategory of unital $\mathcal{V}$-semimonoids is
denoted by $\mathbf{USMonoid}(\mathcal{V}).$
\end{punto}

\begin{punto}
\label{GM}Let $(A,\mu ,\eta )$ be a $\mathcal{V}$-semimonoid. A \emph{right }%
$A$\emph{-semimodule} is a datum $(M,\varrho _{M})$ where $M\in \mathcal{V}$
and $\varrho _{M}:M\bullet A\longrightarrow M$ is a morphism in $\mathcal{V}$
such that the following diagrams are commutative%
\begin{equation*}
\begin{array}{ccc}
\xymatrix{M \bullet A \bullet A \ar[rr]^{\varrho_M \bullet A} \ar[dd]_{M
\bullet \mu } & & M \bullet A \ar[dd]^{\varrho_M}\\ & & \\ M \bullet A
\ar[rr]_{\varrho_M} & & M} &  & \xymatrix{M \bullet A \ar[rr]^{\varrho_M} &
& M \ar[dd]^{\mathbf{\omega}_M}\\ & & & \\ M \bullet \mathbf{I} \ar[uu]^{M
\bullet \eta} & & \mathbf{I} \bullet M \ar[ll]^{\ell_M}}%
\end{array}%
\end{equation*}%
If $M\overset{\mathbf{\mathbf{\omega }}_{M}}{\simeq }\mathbf{I}\bullet M,$
then we say that $M$ is a \emph{unitary} \emph{right} $A$\emph{-semimodule}.
A morphism of right $A$-semimodules is a morphism $f:M\longrightarrow
M^{\prime }$ in $\mathcal{V}$ such that the following diagram is commutative%
\begin{equation*}
\begin{array}{c}
\xymatrix{M \bullet A \ar[rr]^{\varrho_M} \ar[dd]_{f\bullet A} & & M
\ar[dd]^{f }\\ & & \\ M' \bullet A \ar[rr]_{{\varrho}_{M'}} & & M'}%
\end{array}%
\end{equation*}%
The category of right $A$-semimodules is denoted by $\mathbf{S}_{A};$ the
\emph{full} subcategory of unitary right $A$-semimodules is denoted by $%
\mathbf{US}_{A}.$ Analogously, one can define the category $_{A}\mathbf{S}$
of left $A$-semimodules and its \emph{full }subcategory $_{A}\mathbf{US}$ of
unitary left $A$-semimodules.
\end{punto}

\begin{ex}
If $\mathbf{I}\overset{\mathbf{\omega }_{\mathbf{I}}}{\simeq }\mathbf{I}%
\bullet \mathbf{I},$ then $\mathbf{I}$ is a unital $\mathcal{V}$-semimonoid
with%
\begin{equation*}
\mu _{\mathbf{I}}:\mathbf{I}\bullet \mathbf{I}\overset{\mathbf{\omega }_{%
\mathbf{I}}^{-1}}{\longrightarrow }\mathbf{I}\text{ and }\eta _{\mathbf{I}}:%
\mathbf{I}\overset{\mathrm{id}}{\longrightarrow }\mathbf{I}.
\end{equation*}%
Moreover, every $M\in \mathbf{U}\mathcal{V}$ is a unitary $(\mathbf{I},%
\mathbf{I})$-bisemimodule through $\varrho _{M}^{r}:M\bullet \mathbf{I}%
\overset{\rho _{M}}{\simeq }M$ and $\varrho _{M}^{l}:I\bullet M\overset{%
\lambda _{M}}{\simeq }M.$
\end{ex}

\begin{punto}
Let $A$ be a $\mathcal{V}$-semimonoid and $M$ a right $A$-semimodule. We
have a functor%
\begin{equation*}
-\bullet M:\mathcal{V}\longrightarrow \mathbf{S}_{A},
\end{equation*}%
where for any $X\in \mathcal{V}$ we have a structure of a right $A$%
-semimodule on $X\bullet M$ given by%
\begin{equation*}
\varrho _{X\bullet M}:(X\bullet M)\bullet A\overset{{\mathbf{\gamma }}%
_{X,M,A}}{\longrightarrow }X\bullet (M\bullet A)\longrightarrow X\bullet M.
\end{equation*}%
Similarly, if $M$ is a left $A$-semimodule, then we have a functor $M\bullet
-:\mathcal{V}\longrightarrow $ $_{A}\mathbf{S}.$
\end{punto}

\begin{punto}
Let $A$ and $B$ be $\mathcal{V}$-semimonoids. Let $M$ be a left $B$%
-semimodule as well as a right $A$-semimodule and consider $B\bullet M\in
\mathbf{S}_{A}$ and $M\bullet A\in $ $_{B}\mathbf{M}.$ We say that $M$ is a $%
(B,A)$\emph{-bisemimodule} iff $\varrho _{(M;B)}:B\bullet M\longrightarrow M$
is a morphism in $\mathbf{S}_{A},$ equivalently, $\varrho _{(M;A)}:M\bullet
A\longrightarrow M$ is a morphism in $_{B}\mathbf{S}.$ The category of
(unitary) $(A,B)$-bisemimodules with morphisms being in $_{B}\mathbf{S}\cap
\mathbf{S}_{A}$ is denoted by $_{B}\mathbf{S}_{A}$ ($_{B}\mathbf{US}_{A}$).
Indeed, every (unital) $\mathcal{V}$-semimonoid $A$ is a (unitary) $(A,A)$%
-bisemimodule in a canonical way.
\end{punto}

\begin{proposition}
\label{F(A)}Every semiunital semimonoidal functor $F:(\mathcal{V},\bullet ,%
\mathbf{I}_{\mathcal{V}})\longrightarrow (\mathcal{W},\otimes ,\mathbf{I}_{%
\mathcal{W}})$ lifts to a functor%
\begin{equation*}
\widetilde{F}:\mathbf{SMonoid}(\mathcal{V})\longrightarrow \mathbf{SMonoid}(%
\mathcal{W}),\text{ }A\longmapsto F(A)
\end{equation*}%
that commutes with the forgetful functors $U_{\mathcal{V}}:\mathbf{SMon}(%
\mathcal{V})\longrightarrow \mathcal{V}$ and $U_{\mathcal{W}}:\mathbf{SMon}(%
\mathcal{W})\longrightarrow \mathcal{W}.$
\end{proposition}

\begin{Beweis}
Let $(A,\mu _{A},\eta _{A})$ be a semimonoid in $\mathcal{V}$ and consider $%
B:=F(A).$ Define%
\begin{eqnarray*}
\mu _{B} &:&F(A)\otimes F(A)\overset{\phi _{A,A}}{\longrightarrow }%
F(A\bullet A)\overset{F(\mu _{A})}{\longrightarrow }F(A); \\
\eta _{B} &:&\mathbf{I}_{\mathcal{W}}\overset{\phi }{\longrightarrow }F(%
\mathbf{I}_{\mathcal{V}})\overset{F(\eta _{A})}{\longrightarrow }F(A).
\end{eqnarray*}%
One checks easily that $(B,\mu _{B},\eta _{B})$ is a semimonoid in $\mathcal{%
W}.$ If $f:A\longrightarrow A^{\prime }$ is a morphism of $\mathcal{V}$%
-semimonoids, then examining the involved diagrams shows that $%
F(f):F(A)\longrightarrow F(A^{\prime })$ is a morphism of $\mathcal{W}$%
-semimonoids. Finally, it is clear that $U_{\mathcal{W}}\circ \widetilde{F}%
=F\circ U_{\mathcal{V}}.\blacksquare $
\end{Beweis}

\begin{proposition}
\label{monad->alg}Let $(A,\mu ,\eta )$ be a $\mathcal{V}$-semimonoid.

\begin{enumerate}
\item We have $\mathbb{J}$-monads%
\begin{equation*}
-\bullet A:\mathcal{V}\longrightarrow \mathcal{V}\text{ and }A\bullet -:%
\mathcal{V}\longrightarrow \mathcal{V}
\end{equation*}%
and isomorphisms of categories%
\begin{equation*}
\mathbf{S}_{A}\simeq \mathcal{V}_{(-\bullet A;\mathbb{J})}\text{ and }_{A}%
\mathbf{S}\simeq \mathcal{V}_{(-\bullet A;\mathbb{J})}.
\end{equation*}

\item If $B$ is a $\mathcal{V}$-semimonoid, then we have $\mathbb{J}$-monads%
\begin{equation*}
-\bullet A:\text{ }_{B}\mathbf{S}\longrightarrow \text{ }_{B}\mathbf{S}\text{
and }B\bullet -:\mathbf{S}_{A}\longrightarrow \mathbf{S}_{A}
\end{equation*}%
and isomorphisms of categories%
\begin{equation*}
(_{B}\mathbf{S)}_{(-\bullet A;\mathbb{J})}\simeq \text{ }_{B}\mathbf{S}%
_{A}\simeq (\mathbf{S}_{A})_{(B\bullet -;\mathbb{J})}.
\end{equation*}
\end{enumerate}
\end{proposition}

\begin{Beweis}
Consider the natural transformations%
\begin{eqnarray*}
\mathbf{\mu } &:&(-\bullet A)\bullet A\longrightarrow -\bullet A,\text{ }%
\mathbf{\mu }_{X}:(X\bullet A)\bullet A\overset{{\mathbf{\gamma }}_{X,A,A}}{%
\simeq }X\bullet (A\bullet A)\overset{X\bullet \mu }{\longrightarrow }%
X\bullet A, \\
\mathbf{\nu } &:&\mathbb{J}\longrightarrow \mathbf{-}\bullet A,\text{ }%
\mathbf{\nu }_{X}:\mathbf{I}\bullet X\overset{\ell _{X}}{\simeq }X\bullet
\mathbf{I}\overset{X\bullet \eta }{\longrightarrow }X\bullet A.
\end{eqnarray*}%
One can easily check that $(-\bullet A,\mathbf{\mu },\mathbf{\omega },%
\mathbf{\nu })$ is a $\mathbb{J}$-monad. The isomorphism $\mathbf{SM}%
_{A}\simeq \mathcal{V}_{(-\bullet A;\mathbb{J})}$ follows immediately from
comparing the corresponding diagrams. The other assertions can also be
checked easily.$\blacksquare $
\end{Beweis}

An object $G$ in cocomplete category $\mathfrak{A}$ is said to be a (\emph{%
regular})\emph{\ generator }iff for every $X\in \mathfrak{A},$ there exists
a canonical (\emph{regular})\emph{\ epimorphism} $f_{X}:\dbigsqcup\limits_{f%
\in \mathfrak{A}(G,X)}G\longrightarrow X$ \cite[p. 199]{BW2005} (see also
\cite{Kel2005}, \cite{Ver}); recall that an arrow in $\mathfrak{A}$ is said
to be a \emph{regular epimorphism} iff it is a coequalizer (of its \emph{%
kernel pair}).

\begin{theorem}
\label{thm-monad}Let $\mathcal{V}$ be cocomplete, $\mathbf{I}$ and $A\in
\mathcal{V}$ be firm and assume that $\mathbf{I}$ is a regular generator in $%
\mathcal{V}$ and that both $A\bullet -$ and $-\bullet A$ preserve colimits
in $\mathcal{V}.$ There is a bijective correspondence between the semimonoid
structures on $A,$ the $\mathbb{J}$-monad structures $(-\bullet A,\mathbf{%
\mu };\mathbf{\omega },\mathbf{\nu };\mathbb{J})$ and the $\mathbb{J}$-monad
structures $(A\bullet -,\mathbf{\mu },\mathbf{\omega },\mathbf{\nu };\mathbb{%
J}).$
\end{theorem}

\begin{Beweis}
Assume that $(-\bullet A,\mathbf{\mu },\mathbf{\omega },\mathbf{\nu })$ is a
$\mathbb{J}$-monad and consider%
\begin{eqnarray*}
\mu &:&A\bullet A\overset{\mathbf{\omega }_{A}\bullet A}{\longrightarrow }%
\mathbf{I}\bullet A\bullet A\overset{\mathbf{\mu }_{\mathbf{I}}}{%
\longrightarrow }\mathbf{I}\bullet A\overset{\lambda _{A}}{\simeq }A; \\
\eta &:&\mathbf{I}\overset{\mathbf{\omega }_{\mathbf{I}}}{\longrightarrow }%
\mathbf{I}\bullet \mathbf{I}\overset{\mathbf{\nu }_{\mathbf{I}}}{%
\longrightarrow }\mathbf{I}\bullet A\overset{\lambda _{A}}{\simeq }A.
\end{eqnarray*}%
Clearly, $(A,\mu ,\eta )$ is a (unital) semimonoid. The converse follow by
Proposition \ref{monad->alg}. The proof of the bijective correspondence is
similar to that in the proof of \cite[Theorem 3.9]{Ver}. The statement
corresponding to the endo-functor $A\bullet -$ can be proved analogously.$%
\blacksquare $
\end{Beweis}

\subsection*{Semicomonoids}

\begin{punto}
A $\mathcal{V}$\emph{-semicomonoid} is a datum $(C,\Delta ,\varepsilon )$
where $C\in \mathcal{V},$ $\Delta :C\longrightarrow C\bullet C,$ $%
\varepsilon :C\longrightarrow \mathbf{I}$ are morphisms in $\mathcal{V}$
such that the following diagrams are commutative%
\begin{equation*}
\begin{array}{ccc}
\xymatrix{C \ar[rr]^{\Delta} \ar[dd]_{\Delta} & & C \bullet C
\ar[dd]^{\Delta \bullet C}\\ & & \\ C \bullet C \ar[rr]_{C \bullet \Delta} &
& C \bullet C \bullet C} &  & \xymatrix{C \bullet C \ar[dd]_{\varepsilon
\bullet C} & & C \ar[dd]^{{\mathbf{\omega}}_C} \ar[rr]^{\Delta}
\ar[ll]_{\Delta} & & C \bullet C \ar[dd]^{C \bullet \varepsilon}\\ & & & & &
\\ \mathbf{I} \bullet C \ar@{=}[rr] & & \mathbf{I} \bullet C
\ar[rr]_{\ell_C} & & C \bullet \mathbf{I} }%
\end{array}%
\end{equation*}%
If $C\overset{\mathbf{\omega }_{C}}{\simeq }I\bullet C,$ then we say that $C$
is a \emph{counital }$\mathcal{V}$\emph{-semicomonoid}. A \emph{morphism of }%
$\mathcal{V}$\emph{-semicomonoids} $f:(C,\Delta ,\varepsilon
)\longrightarrow (C^{\prime },\Delta ^{\prime },\varepsilon ^{\prime })$ is
a morphism in $\mathcal{V}$ such that the following diagrams are commutative%
\begin{equation*}
\begin{array}{ccc}
\xymatrix{C \ar[rr]^{\Delta} \ar[dd]_{f} & & C \bullet C \ar[dd]^{f \bullet
f}\\ & & \\ C' \ar[rr]_{{\Delta}'} & & C' \bullet C'} &  & \xymatrix{C
\ar[dd]_{f} \ar[rr]^{\varepsilon} & & \mathbf{I} \ar@{=}[dd] & & \\ \\ C'
\ar[rr]_{{\varepsilon}'} & & \mathbf{I} & & }%
\end{array}%
\end{equation*}%
The category of $\mathcal{V}$-semicomonoids is denoted by $\mathbf{SComonoid}%
(\mathcal{V});$ the \emph{full} subcategory of counital $\mathcal{V}$%
-semicomonoids is denoted by $\mathbf{USComonoid}(\mathcal{V}).$
\end{punto}

\begin{punto}
Let $(C,\Delta ,\varepsilon )$ be a $\mathcal{V}$-semicomonoid. A \emph{%
right }$C$\emph{-semicomodule} is a datum $(M,\varrho ^{M})$ where $M\in
\mathcal{V}$ and $\varrho ^{M}:M\longrightarrow M\bullet C$ are morphisms in
$\mathcal{V}$ such that the following diagrams are commutative%
\begin{equation*}
\begin{array}{ccc}
\xymatrix{M \ar[rr]^{\varrho^M } \ar[dd]_{\varrho^M} & & M \bullet C
\ar[dd]^{\varrho^M \bullet C}\\ & & \\ M \bullet C \ar[rr]_{M \bullet
\Delta_C} & & M \bullet C \bullet C} &  & \xymatrix{M
\ar[dd]_{\mathbf{\omega}_M}\ar[rr]^{\varrho^M} & & M \bullet C \ar[dd]^{M
\bullet \varepsilon_C}\\ & & & \\ \mathbf{I} \bullet M & & M \bullet
\mathbf{I} \ar[ll]^{\wp_M}}%
\end{array}%
\end{equation*}%
A \emph{morphism of right }$C$\emph{-semicomodules} is a morphism $%
f:M\longrightarrow M^{\prime }$ in $\mathcal{V}$ such that the following
diagram is commutative%
\begin{equation*}
\begin{array}{c}
\xymatrix{M \ar[rr]^{\varrho^M} \ar[dd]_{f} & & M \bullet C \ar[dd]^{f
\bullet C}\\ & & \\ M' \ar[rr]_{{\varrho}^{M'}} & & M' \bullet C}%
\end{array}%
\end{equation*}%
The category of right $C$-semicomodules is denoted by $\mathbf{S}^{C};$ the
category of \emph{counitary} right $C$-semicomodules is denoted by $\mathbf{%
CS}^{C}.$ Analogously, one can define the category $^{C}\mathbf{S}$ of \emph{%
left }$C$\emph{-semicomodules} and its \emph{full} subcategory $^{C}\mathbf{%
CS}$ of \emph{counitary left }$C$\emph{-semicomodules.}
\end{punto}

\begin{remark}
We prefer to use the terminology unital semimonoids (counital semicomonoids)
over monoids (comonoids) which we reserve for monoidal categories. For
example, the category of unital semimonoids in the monoidal category $%
\mathbf{Set}$ of sets is the category $\mathbf{Monoid}$ of usual monoids of
the sense of Abstract Algebra. The same applies for unitary semimodules
(counitary semicomodules). This is also consistent with the classical
terminology of semirings and semimodules used in Section 5.
\end{remark}

\begin{punto}
Let $C$ be a $\mathcal{V}$-semicomonoid and $M$ a semicounitary right $C$%
-semicomodule. Then we have a functor%
\begin{equation*}
-\bullet M:\mathcal{V}\longrightarrow \mathbf{S}^{C},
\end{equation*}%
where for any $X\in \mathcal{V},$ we have a structure of a right $C$%
-semicomodule on $X\bullet M$ given by
\begin{equation*}
\varrho ^{X\bullet M}:X\bullet M\overset{X\bullet \varrho ^{M}}{%
\longrightarrow }X\bullet (M\bullet C)\overset{{\mathbf{\gamma }}%
_{X,M,C}^{-1}}{\longrightarrow }(X\bullet M)\bullet C.
\end{equation*}%
Similarly, if $M$ is a left $C$-semicomodule, then we have a functor $%
M\bullet -:\mathcal{V}\longrightarrow $ $^{C}\mathbf{S}.$
\end{punto}

\begin{punto}
Let $C$ and $D$ be $\mathcal{V}$-semicomonoids. Let $M$ be a left $D$%
-semicomodule and a right $C$-semicomodule and consider $D\bullet M\in
\mathbf{S}^{C}$ and $M\bullet C\in $ $^{D}\mathbf{S}.$ We say that $M$ is a $%
(D,C)$\emph{-bisemicomodule} iff $\varrho ^{(M;D)}:M\longrightarrow D\bullet
M$ is a morphism in $\mathbf{S}^{C};$ equivalently, $\varrho
^{(M;C)}:M\longrightarrow M\bullet C$ is a morphism in $^{D}\mathbf{S}.$ The
category of $(D,C)$-bisemicomodules with morphisms in $^{D}\mathbf{S}\cap
\mathbf{S}^{C}$ is denoted by $^{D}\mathbf{S}^{C}.$ The \emph{full}
subcategory of counitary $(D,C)$-bisemicomodules is denoted by $^{D}\mathbf{%
CS}^{C}.$ Indeed, every $\mathcal{V}$-semicomonoid $C$ is a $(C,C)$%
-bisemicomodule in a canonical way.
\end{punto}

\begin{ex}
$\mathbf{I}$ is $\mathcal{V}$-semicomonoid with%
\begin{equation*}
\Delta _{\mathbf{I}}:\mathbf{I}\overset{\mathbf{\omega }_{\mathbf{I}}}{%
\longrightarrow }\mathbf{I}\bullet \mathbf{I}\text{ and }\varepsilon _{%
\mathbf{I}}:\mathbf{I}\overset{\mathrm{id}}{\longrightarrow }\mathbf{I}.
\end{equation*}%
Moreover, every (firm) $M\in \mathcal{V}$ is a (counitary) $(\mathbf{I},%
\mathbf{I})$-bisemicomodule with $\varrho _{M}^{l}:M\overset{\mathbf{\omega }%
_{M}}{\longrightarrow }\mathbf{I}\bullet M$ and $\varrho _{M}^{r}:M\overset{%
\ell _{M}\circ \mathbf{\omega }_{M}}{\longrightarrow }M\bullet \mathbf{I}.$
\end{ex}

Dual to Proposition \ref{F(A)}, we have

\begin{proposition}
\label{F(C)}Every semiunital op-semimonoidal functor $F:(\mathcal{V},\bullet
,\mathbf{I}_{\mathcal{V}})\longrightarrow (\mathcal{W},\otimes ,\mathbf{I}_{%
\mathcal{W}})$ lifts to a functor%
\begin{equation*}
\widetilde{F}:\mathbf{SCMonoid}(\mathcal{V})\longrightarrow \mathbf{SCMonoid}%
(\mathcal{W}),\text{ }C\longmapsto F(C)
\end{equation*}%
which commutes with the forgetful functors $U_{\mathcal{V}}:\mathbf{SCMon}(%
\mathcal{V})\longrightarrow \mathcal{V}$ and $U_{\mathcal{W}}:\mathbf{SCMon}(%
\mathcal{W})\longrightarrow \mathcal{W}.$
\end{proposition}

Dual to Proposition \ref{monad->alg}, we obtain

\begin{proposition}
\label{com->coal}Let $(C,\Delta ,\varepsilon )$ be a $\mathcal{V}$%
-semicomonoid.

\begin{enumerate}
\item We have $\mathbb{J}$-comonads%
\begin{equation*}
-\bullet C:\mathcal{V}\longrightarrow \mathcal{V}\text{ and }C\bullet -:%
\mathcal{V}\longrightarrow \mathcal{V}
\end{equation*}%
and isomorphisms of categories%
\begin{equation*}
\mathbf{S}^{C}\simeq \mathcal{V}^{(-\bullet C;\mathbb{J})}\text{ and }^{C}%
\mathbf{S}\simeq \mathcal{V}^{(-\bullet C;\mathbb{J})}.
\end{equation*}

\item If $D$ is a $\mathcal{V}$-semicomonoid, then we have $\mathbb{J}$%
-comonads%
\begin{equation*}
-\bullet C:\text{ }^{D}\mathbf{S}\longrightarrow \text{ }^{D}\mathbf{S}\text{
and }D\bullet -:\mathbf{S}^{C}\longrightarrow \mathbf{S}^{C}
\end{equation*}%
and isomorphisms of categories%
\begin{equation*}
(^{D}\mathbf{M)}^{(-\bullet C;\mathbb{J})}\simeq \text{ }^{D}\mathbf{S}%
^{C}\simeq (\mathbf{S}^{C})^{(D\bullet -;\mathbb{J})}.
\end{equation*}
\end{enumerate}
\end{proposition}

Our second reconstruction result is obtained in a way similar to that of
Theorem \ref{thm-monad}:

\begin{theorem}
\label{thm-comonad}Let $\mathcal{V}$ be cocomplete, $\mathbf{I}$ and $C\in
\mathcal{V}$ be firm and assume that $\mathbf{I}$ is a regular generator and
that both $C\bullet -$ and $-\bullet C$ respect colimits in $\mathcal{V}.$
There is a bijective correspondence between the semicomonoid\ structures on $%
C,$ the $\mathbb{J}$-comonad structures $(-\bullet C,\mathbf{\Delta },%
\mathbf{\omega },\mathbf{\varepsilon };\mathbb{J})$ and the $\mathbb{J}$%
-comonad structures $(C\bullet -,\mathbf{\Delta },\mathbf{\omega },\mathbf{%
\varepsilon };\mathbb{J}).$
\end{theorem}

\begin{proposition}
Let $(C,\Delta ,\varepsilon )$ be a semicomonoid and $(A,\mu ,\eta )$ a
unital semimonoid. Then $(\mathcal{V}(C,A),\ast ,\mathfrak{e})$ is a monoid
in $\mathbf{Set}$ with multiplication and unit given by%
\begin{equation*}
f\ast g:=\mu \circ (f\bullet g)\circ \Delta \text{ and }\mathfrak{e}:=\eta
\circ \varepsilon .
\end{equation*}
\end{proposition}

\begin{proposition}
Let $\varphi :(C,\Delta _{C},\varepsilon _{C})\longrightarrow (D,\Delta
_{D},\varepsilon _{D})$ be a morphism of semicomonoids and $\psi :(A,\mu
_{A},\eta _{A})\longrightarrow (B,\mu _{B},\eta _{B})$ be a morphism of
unital semimonoids. Then%
\begin{equation*}
\mathcal{V}(D,A)\overset{<-,A>}{\longrightarrow }\mathcal{V(}C,A),\text{ }%
f\longmapsto f\circ \varphi \text{ and }\mathcal{V}(C,A)\overset{<C,->}{%
\longrightarrow }\mathcal{V}(C,B),\text{ }g\longmapsto \psi \circ g
\end{equation*}%
are morphisms of monoids in $\mathbf{Set}.$ In particular, we have functors%
\begin{equation*}
\mathcal{V}(-,A):\mathbf{SCMonoid}_{\mathcal{V}}\longrightarrow \mathbf{%
Monoid}\text{ and }\mathcal{V}(C,-):\mathbf{SMonoid}_{\mathcal{V}%
}\longrightarrow \mathbf{Monoid}.
\end{equation*}
\end{proposition}

\section{A concrete example}

\qquad In this section we give applications to the category of bisemimodules
over a base semialgebra.

\subsection*{Semirings and Semimodules}

\qquad For the convenience of the reader and to make the manuscript
self-contained, we begin this section by recalling some \emph{basic}
definitions and results on semirings and their semimodules.

\begin{definition}
A \emph{semiring} is an algebraic structure $(S,+,\cdot ,0,1)$ consisting of
a non-empty set $S$ with two binary operations \textquotedblleft $+$%
\textquotedblright\ (addition) and \textquotedblleft $\cdot $%
\textquotedblright\ (multiplication) satisfying the following axioms:

\begin{enumerate}
\item $(S,+,0)$ is a commutative monoid with neutral element $0_{S};$

\item $(S,\cdot ,1)$ is a monoid with neutral element $1;$

\item $x\cdot (y+z)=x\cdot y+x\cdot z$ and $(y+z)\cdot x=y\cdot x+z\cdot x$
for all $x,y,z\in S;$

\item $0\cdot s=0=s\cdot 0$ for every $s\in S$ (\emph{i.e.} $0$ is \emph{%
absorbing}).
\end{enumerate}
\end{definition}

\begin{punto}
Let $S,S^{\prime }$ be semirings. A map $f:S\longrightarrow S^{\prime }$ is
said to be a \emph{morphism of semirings} iff for all $s_{1},s_{2}\in S:$%
\begin{equation*}
f(s_{1}+s_{2})=f(s_{1})+f(s_{2}),\text{ }f(s_{1}s_{2})=f(s_{1})f(s_{2}),%
\text{ }f(0_{S})=0_{S^{\prime }}\text{ and }f(1_{S})=1_{S^{\prime }}.
\end{equation*}%
The category of semirings is denoted by $\mathbf{SRng}.$
\end{punto}

\begin{punto}
Let $(S,+,\cdot )$ be a semiring. We say that $S$ is

\emph{cancellative} iff the additive semigroup $(S,+)$ is cancellative,
\emph{i.e. }whenever $s,s^{\prime },s^{\prime \prime }\in S$ we have%
\begin{equation*}
s+s^{\prime }=s+s^{\prime \prime }\Rightarrow s^{\prime }=s^{\prime \prime }.
\end{equation*}

\emph{commutative} iff the multiplicative semigroup $(S,\cdot )$ is
commutative;

\emph{semifield} iff $(S\backslash \{0\},\cdot ,1)$ is a commutative group.
\end{punto}

\begin{exs}
Rings are indeed semirings. The first natural example of a (\emph{commutative%
})\emph{\ }semiring which is \emph{not} a ring is $(\mathbb{N}_{0},+,\cdot ),
$ the set of non-negative integers. The semirings $(\mathbb{R}%
_{0}^{+},+,\cdot )$ and $(\mathbb{Q}_{0}^{+},+,\cdot )$ are indeed
semifields. For any associative ring $R$ we have a semiring structure $(%
\mathrm{Ideal}(R),+,\cdot )$ on the set $\mathrm{Ideal}(R)$ of (two-sided)
ideals of $R.$ Any distributive complete lattice $(\mathcal{L},\wedge ,\vee
,0,1)$ is a semiring. For more examples, the reader may refer to \cite%
{Gol1999a}. In the sequel, we always assume that $0_{S}\neq 1_{S}.$
\end{exs}

\begin{definition}
Let $S$ be a semiring. A \emph{right }$S$\emph{-semimodule} is an algebraic
structure $(M,+,0_{M})$ consisting of a non-empty set $M,$ a binary
operation \textquotedblleft $+$\textquotedblright\ along with a right $S$%
-action%
\begin{equation*}
M\times S\longrightarrow M,\text{ }(m,s)\mapsto ms,
\end{equation*}%
such that:

\begin{enumerate}
\item $(M,+,0_{M})$ is a commutative monoid with neutral element $0_{M};$

\item $(ms)s^{\prime }=m(ss^{\prime }),$ $(m+m^{\prime })s=ms+m^{\prime }s$
and $m(s+s^{\prime })=ms+ms^{\prime }$ for all $s,s^{\prime }\in S$ and $%
m,m^{\prime }\in M;$

\item $m1_{S}=m$ and $m0_{S}=0_{M}=0_{M}s$ for all $m\in M$ and $s\in S.$
\end{enumerate}
\end{definition}

\begin{punto}
Let $M,M^{\prime }$ be right $S$-semimodules. A map $f:M\longrightarrow
M^{\prime }$ is said to be a \emph{morphism of }$S$\emph{-semimodules} (or $S
$\emph{-linear}) iff for all $m_{1},m_{2}\in M$ and $s\in S:$%
\begin{equation*}
f(m_{1}+m_{2})=f(m_{1})+f(m_{2})\text{ and }f(ms)=f(m)s.
\end{equation*}%
The set $\mathrm{Hom}_{S}(M,M^{\prime })$ of $S$-linear maps from $M$ to $%
M^{\prime }$ is clearly a commutative monoid under addition. The category of
right $S$-semimodules is denoted by $\mathbb{S}_{S}.$ Analogously, one can
define the category $_{S}\mathbb{S}$ of left $S$-semimodules. A right (left)
$S$-semimodule is said to be \emph{cancellative} iff the semigroup $(M,+)$
is cancellative. With $\mathbb{CS}_{S}\subseteq \mathbb{S}_{S}$ (resp. $_{S}%
\mathbb{CS}\subseteq $ $_{S}\mathbb{S})$ we denote the \emph{full}
subcategory of cancellative right (left) $S$-semimodules. For two semirings $%
S,T,$ an $(S,T)$-bisemimodule $M$ has a structure of a left $S$-semimodule
and a right $T$-semimodule such that $(sm)t=s(mt)$ for all $m\in M,$ $s\in S$
and $t\in T.$ The category of $(S,T)$-bisemimodules and $S$-linear $T$%
-linear maps is denoted by $_{S}\mathbb{S}_{T};$ the \emph{full} subcategory
of cancellative $(S,T)$-bisemimodules is denoted by $_{S}\mathbb{CS}_{T}.$
\end{punto}

\begin{punto}
Let $M$ be a right $S$-semimodule. An $S$-congruence on $M$ is an
equivalence relation $\equiv $ such that%
\begin{equation*}
m_{1}\equiv m_{2}\Rightarrow m_{1}s+m\equiv m_{2}s+m\text{ for all }%
m_{1},m_{2},m\in M\text{ and }s\in S.
\end{equation*}%
In particular, we have an $S$-congruence relation $\equiv _{\lbrack 0]}$ on $%
M$ defined by%
\begin{equation*}
m\equiv _{\lbrack 0]}\text{ }m^{\prime }\text{ }\Longleftrightarrow
m+m^{\prime \prime }=m^{\prime }+m^{\prime \prime }\text{ for some }%
m^{\prime \prime }\in M.
\end{equation*}%
The quotient $S$-semimodule $M/\equiv _{\lbrack 0]}$ is indeed cancellative
and we have a canonical surjection $\mathfrak{c}_{M}:M\longrightarrow
\mathfrak{c}(M),$ where $\mathfrak{c}(M):=M/\equiv _{\lbrack 0]},$ with%
\begin{equation*}
\mathrm{Ker}(\mathfrak{c}_{M})=\{m\in M\mid m+m^{\prime \prime }=m^{\prime
\prime }\text{ for some }m^{\prime \prime }\in M\}.
\end{equation*}%
The class of cancellative right $S$-semimodules is a \emph{reflective}
subcategory of $\mathbb{S}_{S}$ in the sense that the functor $\mathfrak{c}:%
\mathbb{S}_{S}\longrightarrow \mathbb{CS}_{S}$ is left adjoint to the
embedding functor $\mathbb{CS}_{S}\hookrightarrow \mathbb{S}_{S},$ \emph{i.e.%
} for any $S$-semimodule $M$ and any \emph{cancellative} $S$-semimodule $N$
we have a natural isomorphism of commutative monoids $\mathrm{Hom}_{S}(%
\mathfrak{c}(M),N)\simeq \mathrm{Hom}_{S}(M,N)$ \cite[p.517]{Tak1981}.
\end{punto}

\subsubsection*{Takahashi's Tensor-like Product}

\begin{punto}
(\cite[page 187]{Gol1999a}) Let $M_{S}$ be a right $S$-semimodule, $_{S}N$ a
left $S$-semimodule and consider the Abelian monoid $U:=S^{(M\times
N)}\times S^{(M\times N)}.$ Let $U^{\prime }\subseteq S^{(M\times N)}\times
S^{(M\times N)}$ be the \emph{symmetric }$S$-subsemimodule generated by the
set of elements of the form%
\begin{equation*}
\begin{array}{ccc}
(\delta _{(m_{1}+m_{2},n)},\delta _{(m_{1},n)}+\delta _{(m_{2},n)}), &  &
(\delta _{(m_{1},n)}+\delta _{(m_{2},n)},\delta _{(m_{1}+m_{2},n)}), \\
(\delta _{(m,n_{1}+n_{2})},\delta _{(m,n_{1})}+\delta _{(m,n_{2})}), &  &
(\delta _{(m,n_{1})}+\delta _{(m,n_{2})},\delta _{(m,n_{1}+n_{2})}), \\
(\delta _{(ms,n)},\delta _{(m,sn)}), &  & (\delta _{(m,sn)},\delta
_{(ms,n)}),%
\end{array}%
\end{equation*}%
where%
\begin{equation*}
\delta _{m,n}(m,n)=\left\{
\begin{array}{ccc}
1_{S}, &  & m=n \\
&  &  \\
0, &  & m\neq n.%
\end{array}%
\right.
\end{equation*}%
Let $\equiv $ be the $S$-congruence relation on $S^{(M\times N)}$ defined by%
\begin{equation*}
f\equiv f^{\prime }\Longleftrightarrow f+g=f^{\prime }+g^{\prime }\text{ for
some }(g,g^{\prime })\in U^{\prime }.
\end{equation*}%
\emph{Takahashi's tensor-like product} of $M$ and $N$ is defined as $%
M\boxtimes _{S}N:=F/\equiv .$ Notice that there is an $S$\emph{-balanced map}%
\begin{equation*}
\widetilde{\tau }:M\times N\longrightarrow M\boxtimes _{S}N,\text{ }%
(m,n)\mapsto m\boxtimes _{S}n:=(m,n)/\equiv
\end{equation*}%
with the following universal property \cite{Tak1982a}: for every commutative
monoid $G$ and every $S$-bilinear $S$-balanced map $\beta :M\times
N\longrightarrow G$ there exists a \emph{unique} morphism of monoids ${%
\mathbf{\gamma }}:M\boxtimes _{S}N\longrightarrow \mathfrak{c}(G)$ such that
we have a commutative diagram%
\begin{equation}
\xymatrix{M \times N \ar_{\tau }[d] \ar^{\beta}[rr] & & G
\ar[d]^{\mathfrak{c}_G}\\ M \boxtimes_{S} N \ar@{.>}_{{\mathbf{\gamma}}}[rr]
& & \mathfrak{c}(G) }  \label{ten-extend}
\end{equation}
\end{punto}

The following result collects some properties of $-\boxtimes _{S}-$ (compare
with \cite{Abu} and \cite[Proposition 16.15, 16.16]{Gol1999a}):

\begin{proposition}
\label{c-M}Let $M$ be a right $S$-semimodule and $N$ a left $S$-semimodule.

\begin{enumerate}
\item $M\boxtimes _{S}N$ is a \emph{cancellative} commutative monoid.

\item $M_{S}$ \emph{(}$_{S}N$\emph{)} is cancellative if and only if $%
\mathfrak{c}(M)\simeq M$ \emph{(}$\mathfrak{c}(N)\simeq N$\emph{)}. In this
case, we

\item We have natural isomorphisms of functors%
\begin{equation*}
-\boxtimes _{S}S\simeq \mathfrak{c}(-):\text{ }\mathbb{S}_{S}\longrightarrow
\mathbb{S}_{S}\text{ and }S\boxtimes _{S}-\simeq \mathfrak{c}(-):\text{ }_{S}%
\mathbb{S}\longrightarrow \text{ }_{S}\mathbb{S}.
\end{equation*}%
Moreover, we have isomorphisms of functors%
\begin{equation*}
-\boxtimes _{S}S\simeq \mathfrak{c}(-)\simeq S\boxtimes _{S}-:\text{ }_{S}%
\mathbb{S}_{S}\longrightarrow \text{ }_{S}\mathbb{S}_{S}.
\end{equation*}%
We set%
\begin{equation*}
M\boxtimes _{S}S\overset{\vartheta _{M}^{r}}{\simeq }M\text{ and }S\boxtimes
_{S}N\overset{\vartheta _{M}^{l}}{\simeq }N.
\end{equation*}

\item We have idempotent functors%
\begin{equation}
\mathbb{J}:S\boxtimes _{S}-:\text{ }_{S}\mathbb{S}\longrightarrow \text{ }%
_{S}\mathbb{S}\text{ and }\mathbb{K}:=-\boxtimes _{T}T:\mathbb{S}%
_{T}\longrightarrow \mathbb{S}_{T}.  \label{JK}
\end{equation}%
In particular, $\mathfrak{c}(\mathfrak{c}(M))\simeq \mathfrak{c}(M)$ and $%
\mathfrak{c}(\mathfrak{c}(N))\simeq \mathfrak{c}(N).$

\item We have natural isomorphisms of commutative monoids%
\begin{equation}
\mathfrak{c}(M)\boxtimes _{S}N\simeq \mathfrak{c}(M)\boxtimes _{S}\mathfrak{c%
}(N)\simeq M\boxtimes _{S}\mathfrak{c}(N)\simeq M\boxtimes _{S}N\simeq
\mathfrak{c}(M\boxtimes _{S}N).  \label{CJ-l}
\end{equation}
\end{enumerate}
\end{proposition}

\begin{proposition}
\label{adj}Let $S$ and $T$ be semirings, $M$ a right $S$-semimodule and $N$
an $(S,T)$-bisemimodule. Consider the functors%
\begin{equation*}
-\boxtimes _{S}N:\mathbb{S}_{S}\longrightarrow \mathbb{S}_{T},\text{ }%
N\boxtimes _{T}-:\text{ }_{T}\mathbb{S}\longrightarrow \text{ }_{S}\mathbb{S}
\end{equation*}%
and the endo-functors $\mathbb{J}$ and $\mathbb{K}$ in \emph{(\ref{JK})}.

\begin{enumerate}
\item $(-\boxtimes _{S}N,\mathrm{Hom}_{-T}(N,-))$ is a $(\mathbb{J},\mathbb{K%
})$-adjoint pair.

\item $(N\boxtimes _{T}-,\mathrm{Hom}_{S-}(N,-))$ is a $(\mathbb{K},\mathbb{J%
})$-adjoint pair.
\end{enumerate}
\end{proposition}

\begin{Beweis}
For every right\emph{\ }$T$-semimodule $G$ we have natural isomorphisms%
\begin{equation*}
\begin{tabular}{lllll}
$\mathrm{Hom}_{-T}(\mathbb{J}(M)\boxtimes _{S}N),\mathbb{K}(G))$ & $\simeq $
& $\mathrm{Hom}_{-T}(\mathfrak{c}(M)\boxtimes _{S}N),\mathfrak{c}(G))$ &  &
\\
& $\simeq $ & $\mathrm{Hom}_{-T}(M\boxtimes _{S}N,\mathfrak{c}(G))$ &  &  \\
& $\simeq $ & $\mathrm{Hom}_{-S}(M,\mathrm{Hom}_{-T}(N,\mathfrak{c}(G)))$ & (%
\cite[16.15]{Gol1999a}) &  \\
& $\simeq $ & $\mathrm{Hom}_{-S}(\mathfrak{c}(M),\mathrm{Hom}_{-T}(N,%
\mathfrak{c}(G)))$ & (\cite[p. 517]{Tak1981}) &  \\
& $\simeq $ & $\mathrm{Hom}_{-S}(\mathbb{J}(M),\mathrm{Hom}_{-T}(N,\mathbb{K}%
(G))).$ &  &
\end{tabular}%
\end{equation*}%
The second statement can be proved symmetrically.$\blacksquare $
\end{Beweis}

\subsection*{Semiunital Semirings and Semicounitary Semimodules}

\qquad In what follows, $S$ denotes a \emph{commutative} semiring with $%
1_{S}\neq 0_{S},$ $A$ is an $S$\emph{-semialgebra} (\emph{i.e. }a semiring
with a morphism of semirings $\iota _{A}:S\longrightarrow A$), $_{A}\mathbb{S%
}_{A}$ is the category of $(A,A)$-bisemimodules and $_{A}\mathbb{CS}_{A}$ is
its \emph{full} subcategory of cancellative $(A,A)$-bisemimodules. Moreover,
we fix the idempotent endo-functor%
\begin{equation*}
\mathbb{J}:=\mathfrak{c}(-)\simeq A\boxtimes _{A}-\simeq -\boxtimes _{A}A:%
\text{ }_{A}\mathbb{S}_{A}\longrightarrow \text{ }_{A}\mathbb{S}_{A}.
\end{equation*}

Summarizing the observations above, we obtain

\begin{theorem}

\begin{enumerate}
\item $(_{A}\mathbb{S}_{A},\boxtimes ,A)$ is a closed semiunital
semimonoidal category.

\item $(_{A}\mathbb{C}\mathbb{S}_{A},\boxtimes ,\mathfrak{c}(A))$ is a
closed monoidal category.
\end{enumerate}
\end{theorem}

\begin{punto}
By{\normalsize \ }a \emph{semiunital }$A$\emph{-semiring}\ we mean an $(A,A)$%
-bisemimodule $\mathcal{A}$ associated with $(A,A)$-bilinear maps $\mu
\mathcal{_{A}}:\mathcal{A}\boxtimes _{A}\mathcal{A}\longrightarrow \mathcal{A%
}$ and $\eta _{\mathcal{A}}:A\longrightarrow \mathcal{A}$ such that the
following diagrams are commutative%
\begin{equation*}
\begin{tabular}{lll}
$\xymatrix{{\mathcal{A}} \boxtimes_A {\mathcal{A}} \boxtimes_A {\mathcal{A}}
\ar^{\mu_{{\mathcal{A}}} \boxtimes_A \mathcal{A}}[rr] \ar_{\mathcal{A}
\boxtimes_A \mu_{{\mathcal{A}}}}[d] & & {\mathcal{A}} \boxtimes_A
{\mathcal{A}} \ar^{\mu_{{\mathcal{A}}}}[d]\\ {\mathcal{A}} \boxtimes_A
{\mathcal{A}} \ar_{\mu_{{\mathcal{A}}}}[rr] & & {\mathcal{A}}}$ &  & $%
\xymatrix{{\mathcal {A}} \boxtimes_{A} {\mathcal {A}} \ar^{{\mu_{\mathcal
{A}}}}[r] & {\mathcal{A}} \ar_{\mathfrak{c}_{A}}[d] & {\mathcal {A}}
\boxtimes_{A} {\mathcal {A}} \ar_{{\mu _{\mathcal {A}}}}[l] \\ A
\boxtimes_{A} {\mathcal {A}} \ar^{\eta_{\mathcal {A}} \boxtimes_A
{\mathcal{A}}}[u] \ar[r]_(.45){\vartheta_{\mathcal {A}} ^l} &
\mathfrak{c}({\mathcal {A}}) & \mathcal{A} \boxtimes_{A} A \ar[u]_{{\mathcal
{A}} \boxtimes_A \eta_{\mathcal {A}}} \ar[l]^(.45){\vartheta _{\mathcal {A}}
^r} }$%
\end{tabular}%
\end{equation*}%
Let $\mathcal{A}$ and $\mathcal{A}^{\prime }$ be semiunital $A$-semirings.
An $(A,A)$-bilinear map $f:\mathcal{A}\longrightarrow \mathcal{A}^{\prime }$
is called a \emph{morphism of semiunital }$A$-\emph{semirings} iff%
\begin{equation*}
f\circ \mu _{\mathcal{A}}=\mu _{\mathcal{A}^{\prime }}\circ (f\boxtimes
_{A}f)\text{ and }f\circ \eta _{\mathcal{A}}=\eta _{\mathcal{A}^{\prime }}.
\end{equation*}%
The set of morphisms of semiunital $A$-semirings form $\mathcal{A}$ to $%
\mathcal{A}^{\prime }$ is denoted by $\mathrm{SSRng}_{A}(\mathcal{A},%
\mathcal{A}^{\prime }).$ The category of semiunital $A$-semirings will be
denoted by $\mathbf{SSRng}_{A}.$ Indeed, we have an isomorphism of
categories $\mathbf{SSRng}_{A}\simeq \mathbf{SMonoid}(_{A}\mathbb{S}_{A}).$
\end{punto}

\begin{punto}
Let $\mathcal{A}$ be a semiunital $A$-semiring. A \emph{semiunitary right }$%
\mathcal{A}$\emph{-semimodule} is a right $A$-semimodule along with a right $%
A$-linear map $\varrho _{M}:M\boxtimes _{A}\mathcal{A}\longrightarrow M$
such that the following diagrams are commutative%
\begin{equation*}
\begin{array}{ccc}
\xymatrix{M \boxtimes_A \mathcal{A} \boxtimes_A \mathcal{A}
\ar[rr]^{\varrho_M \boxtimes_A \mathcal{A}} \ar[dd]_{M \boxtimes_A
\mu_{\mathcal{A}} } & & M \boxtimes_A \mathcal{A} \ar[dd]^{\varrho_M}\\ & &
\\ M \boxtimes_A \mathcal{A} \ar[rr]_{\varrho_M} & & M} &  & \xymatrix{M
\boxtimes_A A \ar[dd]_{\vartheta_M ^r} \ar[rr]^{M \boxtimes_A
\eta_{\mathcal{A}}} & & M \boxtimes_A \mathcal{A} \ar[dd]^{\varrho_M}\\ & &
& \\ \mathfrak{c}(M) & & M \ar[ll]^{\mathfrak{c}_M}}%
\end{array}%
\end{equation*}%
A \emph{morphism of semiunitary right }$\mathcal{A}$\emph{-semimodules} ($%
\mathcal{A}$-linear) is an $A$-linear map $f:M\longrightarrow M^{\prime }$
such that the following diagram is commutative%
\begin{equation*}
\begin{array}{c}
\xymatrix{M \boxtimes_A \mathcal{A} \ar[rr]^{\varrho_M}
\ar[dd]_{f\boxtimes_A \mathcal{A}} & & M \ar[dd]^{f }\\ & & \\ M'
\boxtimes_A \mathcal{A} \ar[rr]_{{\varrho}_{M'}} & & M'}%
\end{array}%
\end{equation*}%
The category of semiunitary right $\mathcal{A}$-semimodules and $\mathcal{A}$%
-linear maps is denoted by $\mathbb{SS}_{\mathcal{A}}.$ Analogously, one can
define the category $_{\mathcal{A}}\mathbb{S}\mathbb{S}$ of \emph{semiunital
left }$\mathcal{A}$\emph{-semimodules}. For two semiunital $A$-semirings $%
\mathcal{A}$ and $\mathcal{B},$ one can define the category $_{\mathcal{B}}%
\mathbb{S}\mathbb{S}_{\mathcal{A}}$ of $(\mathcal{B},\mathcal{A})$%
-bisemimodules in the obvious way. Considering semiunital $A$-semirings as
semimonoids in $_{A}\mathbb{S}_{A},$ we have indeed isomorphisms of
categories%
\begin{equation}
\mathbb{SS}_{\mathcal{A}}\simeq \mathbf{S}_{\mathcal{A}},\text{ }_{\mathcal{B%
}}\mathbb{SS}\simeq \text{ }_{\mathcal{B}}\mathbf{S},\text{ }_{\mathcal{B}}%
\mathbb{SS}_{\mathcal{A}}\simeq \text{ }_{\mathcal{B}}\mathbf{S}_{\mathcal{A}%
},\text{ }\mathbb{CS}_{\mathcal{A}}\simeq \mathbf{US}_{\mathcal{A}},\text{ }%
_{\mathcal{B}}\mathbb{CS}\simeq \text{ }_{\mathcal{B}}\mathbf{US},\text{ }_{%
\mathcal{B}}\mathbb{CS}_{\mathcal{A}}\simeq \text{ }_{\mathcal{B}}\mathbf{US}%
_{\mathcal{A}}.
\end{equation}
\end{punto}

\begin{remark}
We use semiunital $A$-semirings to stress that such semimonoids are defined
in the semiunital semimonoidal category $(_{A}\mathbb{S}_{A},\boxtimes
_{A},A)$ and to avoid confusion with (unital) $A$\emph{-semirings} which can
be defined as monoids in the monoidal category $(_{A}\mathbb{S}_{A},\otimes
_{A},A).$ The same applies for semicounitary $A$-semicorings below.
\end{remark}

\begin{punto}
Being a variety, in the sense of Universal Algebra, the category $_{A}%
\mathbb{S}_{A}$ of $(A,A)$-bisemimodules is cocomplete. The class of regular
epimorphism in $_{A}\mathbb{S}_{A}$ coincides with that of surjective $(A,A)$%
-bilinear maps. For every $(A,A)$-bisemimodule $M,$ there is a surjective $%
(A,A)$-bilinear map from a free $(A,A)$-bisemimodule to $M$ (compare with
\cite[Proposition 17.11]{Gol1999a}); whence, $A$ is a regular generator.
Moreover, for any $(A,A)$-bisemimodule $X,$ both $X\boxtimes _{A}-,$ $%
-\boxtimes _{A}X:$ $_{A}\mathbb{S}_{A}\longrightarrow $ $_{A}\mathbb{CS}_{A}$
respect colimits since they are left adjoints \cite[Corollary 4.5]{Tak1982a}.
\end{punto}

Applying Theorem \ref{thm-monad} to $_{A}\mathbb{S}_{A},$ we obtain:

\begin{corollary}
\label{thm-A-cancl}Let $A$ be cancellative and $\mathcal{A}$ a cancellative $%
(A,A)$-bisemimodule. There is a bijective correspondence between the
structures of unital $A$-semirings on $\mathcal{A},$ $\mathfrak{c}$-monads
on $\mathcal{A}\boxtimes _{A}-$ and $\mathfrak{c}$-monads $-\boxtimes _{A}%
\mathcal{A}.$
\end{corollary}

\subsection*{Semicounital Semicorings and Semicounitary Semicomodules}

\begin{punto}
A \emph{semicounital }$A$\emph{-semicoring} is an $(A,A)$-bisemimodule
associated with $(A,A)$-bilinear maps $\Delta _{\mathcal{C}}:\mathcal{C}%
\longrightarrow \mathcal{C}\boxtimes _{A}\mathcal{C}$ and $\varepsilon _{%
\mathcal{C}}:\mathcal{C}\longrightarrow A$ such that the following diagrams
are commutative%
\begin{equation}
\begin{tabular}{lll}
$\xymatrix{{\mathcal {C}} \ar^(.45){\Delta_{\mathcal C}}[rr]
\ar_(.45){\Delta_{\mathcal C}}[d] & & {\mathcal {C}} \boxtimes_A {\mathcal
{C}} \ar^(.45){\mathcal {C} \boxtimes_A \Delta_{\mathcal C}}[d]\\ {\mathcal
{C}} \boxtimes_A {\mathcal {C}} \ar_(.45){\Delta_{\mathcal C} \boxtimes_A
\mathcal {C}}[rr] & & {\mathcal {C}} \boxtimes_A {\mathcal {C}} \boxtimes_A
{\mathcal {C}} }$ &  & $\xymatrix{{\mathcal {C}} \boxtimes_A {\mathcal {C}}
\ar_{\varepsilon _{\mathcal {C}} \boxtimes_A \mathcal {C}}[d] & {\mathcal
{C}} \ar_(.4){{\Delta _{\mathcal {C}}}}[l] \ar^(.4){{\Delta _{\mathcal
{C}}}}[r] \ar_{\mathfrak{c}_{C}}[d] & {\mathcal {C}} \boxtimes_A {\mathcal
{C}} \ar[d]^{\mathcal {C} \boxtimes_A \varepsilon _{\mathcal {C}}} \\ A
\boxtimes_A {\mathcal {C}} \ar[r]_(.45){\vartheta_{\mathcal {C}} ^l} &
\mathfrak{c}({\mathcal {C}}) & {\mathcal {C}} \boxtimes_A {A}
\ar[l]^(.45){\vartheta _{\mathcal {C}} ^r} }$%
\end{tabular}
\label{Diag-WC}
\end{equation}%
The map $\Delta _{\mathcal{C}}$ ($\varepsilon _{\mathcal{C}})$ is called the
\emph{comultiplication} (\emph{counity}) of $\mathcal{C}.$ Using
Sweedler-Heyneman's notation, we have for every $c\in \mathcal{C}:$%
\begin{eqnarray*}
\sum c_{11}\boxtimes _{A}c_{12}\boxtimes _{A}c_{2} &=&\sum c_{1}\boxtimes
_{A}c_{21}\boxtimes _{A}c_{22}; \\
\mathfrak{c}(\sum c_{1}\varepsilon _{\mathcal{C}}(c_{2})) &=&\mathfrak{c}%
_{M}(c)=\mathfrak{c}(\sum \varepsilon _{\mathcal{C}}(c_{1})c_{2}).
\end{eqnarray*}%
Let $(\mathcal{C},\Delta ,\varepsilon )$ and $(\mathcal{C}^{\prime },\Delta
^{\prime },\varepsilon ^{\prime })$ be semicounital $A$-semicorings. We call
an $(A,A)$-bilinear map $f:\mathcal{C}\longrightarrow \mathcal{C}^{\prime }$%
\ a \emph{morphism of }$A$\emph{-semicorings} iff%
\begin{equation*}
(f\boxtimes _{A}f)\circ \Delta _{\mathcal{C}}=\Delta _{\mathcal{C}^{\prime
}}\circ f\text{ and }\varepsilon _{\mathcal{C}^{\prime }}\circ f=\varepsilon
_{\mathcal{C}}.
\end{equation*}%
The set of $A$-semicoring morphisms from $\mathcal{C}$ to $\mathcal{C}%
^{\prime }$ is denoted by $\mathrm{SSCog}_{A}(\mathcal{C},\mathcal{C}%
^{\prime }).$ The category of semicounital $A$-semicorings is denoted by $%
\mathbf{SSCrng}_{A}.$ Indeed, we have an isomorphism of categories $\mathbf{%
SSCrng}_{A}\simeq \mathbf{SCMonoid}(_{A}\mathbb{S}_{A}).$
\end{punto}

\begin{punto}
Let $(\mathcal{C},\Delta ,\varepsilon )$ be an $A$-semicoring. A \emph{%
semicounitary right }$\mathcal{C}$\emph{-semicomodule} is a right $A$%
-semimodule $M$ associated with an $A$-linear map%
\begin{equation*}
\varrho ^{M}:M\longrightarrow M\boxtimes _{A}\mathcal{C},\text{ }m\mapsto
\sum m_{<0>}\boxtimes _{A}m_{<1>},
\end{equation*}%
such that the following diagrams are commutative%
\begin{equation*}
\begin{tabular}{lll}
$\xymatrix{M \ar^(.4){\varrho _M}[rr] \ar_(.45){\varrho _M}[d] & & M
\boxtimes_{A} {\mathcal{C}} \ar^(.45){M \boxtimes_A \Delta_{\mathcal{C}}}[d]
\\ M \boxtimes_{A} {\mathcal{C}} \ar_(.4){\varrho _M \boxtimes_A
\mathcal{C}}[rr] & & M \boxtimes_{A} {\mathcal {C}} \boxtimes_{A}
{\mathcal{C}}}$ &  & $\xymatrix{M \ar^(.4){\varrho _M}[rr]
\ar_(.45){\mathfrak{c}_M}[d] & & M \boxtimes_{A} {\mathcal{C}} \ar^(.45){M
\boxtimes_A \varepsilon_{\mathcal{C}}}[d] \\ \mathfrak{c}(M) & & M
\boxtimes_{A} A \ar^(.4){\vartheta _{M}^{r}}[ll] }$%
\end{tabular}%
\end{equation*}%
Using Sweedler-Heyneman's notation, we have for every $m\in M:$%
\begin{eqnarray*}
\sum m_{<0>}\boxtimes _{A}m_{<1>1}\boxtimes _{A}m_{<1>2} &=&\sum
m_{<0><0>}\boxtimes _{A}m_{<0><1>}\boxtimes _{A}m_{<1>}; \\
\mathfrak{c}(\sum m_{<0>}\varepsilon _{\mathcal{C}}(m_{<1>})) &=&\mathfrak{c}%
_{M}(m).
\end{eqnarray*}%
For semicounitary right $\mathcal{C}$-comodules $M,M^{\prime },${\normalsize %
\ }we call an $A$-linear map $f:M\longrightarrow M^{\prime }$ a \emph{%
morphism of semicounitary right }$\mathcal{C}$\emph{-semicomodules}%
{\normalsize \ }(or $\mathcal{C}$\emph{-colinear}) iff the following diagram
is commutative%
\begin{equation*}
\xymatrix{M \ar[rr]^{f} \ar[d]_{\varrho_M} & & N \ar[d]^{\varrho_N}\\ M
\boxtimes_{A} {\mathcal{C}} \ar[rr]_{f \boxtimes_A \mathcal{C}} & & N
\boxtimes_{A} {\mathcal{C}} }
\end{equation*}%
The category of semicounitary right $\mathcal{C}$-semicomodules and $%
\mathcal{C}$-colinear maps is denoted by $\mathbb{SS}^{\mathcal{C}};$ the
\emph{full} subcategory of counitary right $\mathcal{C}$-semicomodules is
denoted by $\mathbb{CS}^{\mathcal{C}}.$ Analogously, one can define the
category $^{\mathcal{C}}\mathbb{SS}$ of \emph{semicounitary left }$\mathcal{C%
}$\emph{-semicomodules} and its \emph{full} subcategory of \emph{counitary
left }$\mathcal{C}$\emph{-semicomodules}. For two semicounital $A$%
-semicorings $\mathcal{C}$ and $\mathcal{D}$ one can define the category $^{%
\mathcal{D}}\mathbb{SS}^{\mathcal{C}}$ of \emph{semicounitary }$(\mathcal{D},%
\mathcal{C})$\emph{-bisemicomodules} and its \emph{full} subcategory of
\emph{counitary }$(\mathcal{D},\mathcal{C})$\emph{-bisemicomodules} in the
obvious way. Considering semicounital $A$-semicorings as semicomonoids in $%
_{A}\mathbb{S}_{A},$ we have indeed isomorphisms of categories%
\begin{equation}
\mathbb{SS}^{\mathcal{C}}\simeq \mathbf{S}^{\mathcal{C}},^{\mathcal{D}}%
\mathbb{SS}\simeq \text{ }^{\mathcal{D}}\mathbf{S},\text{ }^{\mathcal{D}}%
\mathbb{SS}^{\mathcal{C}}\simeq \text{ }^{\mathcal{D}}\mathbf{S}^{\mathcal{C}%
},\text{ }\mathbb{CS}^{\mathcal{C}}\simeq \mathbf{CS}^{\mathcal{C}},^{%
\mathcal{D}}\mathbb{CS}\simeq \text{ }^{\mathcal{D}}\mathbf{CS},\text{ }^{%
\mathcal{D}}\mathbb{CS}^{\mathcal{C}}\simeq \text{ }^{\mathcal{D}}\mathbf{CS}%
^{\mathcal{C}}.
\end{equation}
\end{punto}

Applying Theorem \ref{thm-comonad} to $_{A}\mathbb{S}_{A},$ we obtain:

\begin{corollary}
\label{thm-C-cancel}Let $A$ be cancellative and $\mathcal{C}$ a cancellative
$(A,A)$-bisemimodule. There is a bijective correspondence between the
structures of counital $A$-semicorings on $\mathcal{C},$ $\mathfrak{c}$%
-comonads on $\mathcal{C}\boxtimes _{A}-$ and $\mathfrak{c}$-comonads on $%
-\boxtimes _{A}\mathcal{C}.$
\end{corollary}

Almost all structures of corings over rings (\emph{e.g.} \cite{Abu}, \cite%
{BW2003}) can be transferred to obtain structures of semicorings over
semirings.

\begin{ex}
Let $\kappa :B\longrightarrow A$ be an extension of $S$-semialgebras and
consider $B$ as a $(B,B)$-bisemimodule in the canonical way. One can define
Sweedler's counital $A$-semicoring $\mathcal{C}:=(A\boxtimes _{B}A,\Delta
,\varepsilon )$ with%
\begin{eqnarray*}
\Delta  &:&A\boxtimes _{B}A\longrightarrow (A\boxtimes _{B}A)\boxtimes
_{A}(A\boxtimes _{B}A),\text{ }a\boxtimes _{B}\widetilde{a}\mapsto
(a\boxtimes _{B}1_{A})\boxtimes _{A}(1_{A}\boxtimes _{B}\widetilde{a}); \\
\varepsilon  &:&A\boxtimes _{B}A\longrightarrow A,\text{ }a\boxtimes _{B}%
\widetilde{a}\mapsto a\widetilde{a}.
\end{eqnarray*}
\end{ex}

\textbf{Acknowledgement: }The author thanks the anonymous referee for
her/his careful reading of the paper and, in particular, for suggestions to
fix some results in the initial version of the paper. He also thanks Prof.
Robert Wisbauer for the fruitful discussions on monads and comonads during
his visit to the University of D\"{u}sseldorf in Summer 2009. He also thanks
\emph{Deutsche Akademische Austausch Dienst} (DAAD) for supporting that
visit.

\end{document}